\documentclass[10pt]{article}
\usepackage{amsfonts}
\usepackage[leqno]{amsmath}
\usepackage{graphicx}
\usepackage{latexsym}
\usepackage{amsmath,amsfonts,amssymb,amsthm,mathrsfs,euscript,makeidx,color}
\usepackage{enumerate}
\allowdisplaybreaks[1]
\usepackage{multirow}

\def\sqr#1#2{{\vcenter{\vbox{\hrule height.#2pt
				\hbox{\vrule width.#2pt height#1pt \kern#1pt \vrule width.#2pt}
				\hrule height.#2pt}}}}
\def\signed #1{{\unskip\nobreak\hfil\penalty50
		\hskip2em\hbox{}\nobreak\hfil#1
		\parfillskip=0pt \finalhyphendemerits=0 \par}}
\def\endpf{\signed {$\sqr69$}}

\def\3n{\negthinspace \negthinspace \negthinspace }
\def\2n{\negthinspace \negthinspace }
\def\1n{\negthinspace }
\def\bel{\begin{equation}\label}

\def\dbE{\mathbb{E}}
\def\dbF{\mathbb{F}}

\def\dbP{\mathbb{P}}

\def\dbR{\mathbb{R}}
\def\dbS{\mathbb{S}}

\def\sU{\mathscr{U}}

\def\ds{\displaystyle}

\def\ns{\noalign{\ss}}
\def\e{\varepsilon}
\def\z{\zeta}

\def\l{\lambda}
\def\m{\mu}
\def\n{\nu}
\def\si{\sigma}
\def\t{\tau}
\def\f{\varphi}

\def\i{\infty}
\def\G{\Gamma}

\def\Th{\Theta}
\def\L{\Lambda}

\def\F{\Phi}
\def\O{\Omega}
\def\cF{{\cal F}}

\def\cN{{\cal N}}

\def\cX{{\cal X}}

\def\G{\Gamma}

\def\Th{\Theta}
\def\L{\Lambda}

\def\F{\varPhi}
\def\O{\Omega}

\def\no{\noindent}

\def\ss{\smallskip}
\def\ms{\medskip}

\def\q{\quad}
\def\qq{\qquad}
\def\hb{\hbox}

\def\h1{\outline{$1$}}

\def\hTh{\outline{$\Theta$}}

\def\hh1{\outline{$1$}}
\def\hh2{\outline{$2$}}
\def\hh3{\outline{$3$}}
\def\hh4{\outline{$4$}}
\def\hh5{\outline{$5$}}
\def\hh6{\outline{$6$}}
\def\hh7{\outline{$7$}}
\def\hh8{\outline{$8$}}
\def\hh9{\outline{$9$}}
\def\hh0{\outline{$0$}}

\def\Ra{\mathop{\Rightarrow}}

\def\lan{{\langle}}
\def\ran{{\rangle}}

\def\h{\widehat}
\def\wt{\widetilde}

\def\cd{\cdot}

\def\as{\hbox{\rm a.s.}}

\def\les{\leqslant}
\def\ges{\geqslant}

\def\({\Big (}
\def\){\Big )}
\def\[{\Big[}
\def\]{\Big]}
\def\lan{\langle}
\def\ran{\rangle}

\def\bde{\begin{definition}\label}
	\def\ede{\end{definition}}
\def\bel{\begin{equation}\label}
		\def\ee{\end{equation}}
	\def\bt{\begin{theorem}\label}
		\def\et{\end{theorem}}
	\def\bc{\begin{corollary}\label}
		\def\ec{\end{corollary}}
	\def\bl{\begin{lemma}\label}
		\def\el{\end{lemma}}
	\def\bp{\begin{proposition}\label}
		\def\ep{\end{proposition}}
	\def\bex{\begin{example}\label}
		\def\ex{\end{example}}
	\def\bas{\begin{assumption}}
		\def\eas{\end{assumption}}
	\def\br{\begin{remark}\label}
		\def\er{\end{remark}}
	\def\ba{\begin{array}}
		\def\ea{\end{array}}
	\def\ed{\end{document}}

\def\rf{\eqref}

\def\square#1{\vbox{\hrule\hbox{\vrule height#1%
			\kern#1\vrule}\hrule}}
\def\rectangle#1#2{\vbox{\hrule\hbox{\vrule height#1%
			\kern#2\vrule}\hrule}}

\font\tenbb=msbm10 \font\sevenbb=msbm7 \font\fivebb=msbm5

\newfam\bbfam
\scriptscriptfont\bbfam=\fivebb \textfont\bbfam=\tenbb
\scriptfont\bbfam=\sevenbb

\newtheorem{theorem}{Theorem}[section]
\newtheorem{corollary}[theorem]{Corollary}

\newtheorem{lemma}[theorem]{Lemma}
\newtheorem{proposition}[theorem]{Proposition}

\theoremstyle{definition}
\newtheorem{definition}[theorem]{Definition}

\newtheorem{remark}[theorem]{Remark}

\newtheorem{example}[theorem]{Example}

\makeatletter

\@addtoreset{equation}{section}
\makeatother

\usepackage{xcolor}
\makeatletter

\input pdf-trans
\newbox\qbox
\def\usecolor#1{\csname\string\color@#1\endcsname\space}
\newcommand\bordercolor[1]{\colsplit{1}{#1}}
\newcommand\fillcolor[1]{\colsplit{0}{#1}}
\newcommand\outline[1]{\leavevmode%
	\def\maltext{#1}%
	\setbox\qbox=\hbox{\maltext}%
	\boxgs{Q q 2 Tr \thickness\space w \fillcol\space \bordercol\space}{}%
	\copy\qbox%
}
\makeatother
\newcommand\colsplit[2]{\colorlet{tmpcolor}{#2}\edef\tmp{\usecolor{tmpcolor}}%
	\def\tmpB{}\expandafter\colsplithelp\tmp\relax%
	\ifnum0=#1\relax\edef\fillcol{\tmpB}\else\edef\bordercol{\tmpC}\fi}
\def\colsplithelp#1#2 #3\relax{%
	\edef\tmpB{\tmpB#1#2 }%
	\ifnum `#1>`9\relax\def\tmpC{#3}\else\colsplithelp#3\relax\fi
}
\bordercolor{black}
\fillcolor{white}
\def\thickness{.3}

\def\hTh{\outline{$\Theta$}}

\begin{document}

\title{\bf Stochastic Optimal Linear Quadratic Controls with\\ A Recursive Cost Functional in Infinite Horizon}

\author{Lin Li\thanks{Department of Mathematics, University of Central Florida, Orlando, FL 32816, USA (Email: {\tt lin.li@ucf.edu}).
                           }
~~~and~~~
Jiongmin Yong\thanks{Department of Mathematics, University of Central Florida, Orlando, FL 32816, USA
                    (Email: {\tt jiongmin.yong@ucf.edu}).
                    This author is supported by NSF grant DMS-2305475.}
}

\date{}
\maketitle
\ms

\no{\bf Abstract.} This paper is concerned with a stochastic linear quadratic (LQ, for short) control problem with a recursive cost functional in an infinite horizon. A main difficult is well-posedness of the BSDE in $L^1$ and in infinite horizon. A notion of weighted $L^2$-stabilizability is introduced and characterized, which will lead to an equivalence of the optimal control problem having recursive cost functional with a classical LQ problem. Then all the results of classical problems for open-loop and closed-loop solvability of such an LQ problem can be translated, in terms of the problem under consideration. Then they are characterized as the solvability of a forward-backward stochastic differential equation (FBSDE, for short) and that of algebraic Riccati equation ARE for short). Finally, the nonhomogeneous is briefly discussed.

\ms

\no{\bf Keywords.} Optimal linear quadratic control, recursive cost functional, forward-backward stochastic differential equation, Riccati equation, backward differential equation in $L^1$ space.

\ms

\no{\bf AMS 2020 Mathematics Subject Classification.} 93E20, 49N10, 60H10.

\section{Introduction.}\label{Sec:Intro}

Let $(\O,\cF,\dbP)$ be a complete probability space, with $\cN$ being the set of all $\dbP$-null sets (and their subsets), on which a one-dimensional standard Brownian motion is defined (the case of multi-dimensional standard Brownian motion will be similar). Define
\bel{si}\cF^t_s=\si\(\{W(\t')-W(\t)\bigm|t\les\t<\t'\les s\}\vee\cN\),\qq\cF_t^t=\si(\cN),\qq\cF_\i^t=\bigcup_{s\ges t}\cF_s^t,\ee
where $\si(M)$ (traditionally) stands for the $\si$-field generated by the set $M$. Denote the corresponding filtration by $\dbF^t=\{\cF^t_s\}_{s\ges t}$. In what follows, we will denote $\dbE_s^t[\,\cd\,]=\dbE[\,\cd\,|\cF^t_s]$ (and thus $\dbE^t_t[\,\cd\,]=\dbE[\,\cd\,|\cF^t_t]=\dbE[\,\cd\,]$). When $t=0$, we simply write $\dbF$ and $\cF_T$ instead of $\dbF^0$ and $\cF^0_T$. Now, in $(\O,\cF,\dbF^t,\dbP)$, we begin with the following controlled linear autonomous homogeneous stochastic differential equation (SDE, for short), called {\it state equation}, on any finite horizon $[t,T]$:
\bel{state}\left\{\2n\ba{ll}
\ds dX(s)=[AX(s)+Bu(s)]ds+[CX(s)+Du(s))]dW(s),\q s\in[t,T],\\
\ns\ds X(t)=x.\ea\right.\ee
In the above, $(t,x)\in[0,T]\times\dbR^n$ is called an {\it initial pair}, $X(\cd)$, $u(\cd)$ and $(X(\cd),u(\cd))$ are called the {\it state process}, the {\it control process} and the {\it state-control pair}, respectively. Such a system is denoted by $[A,C;B,D]$ (We are going to use the same notation for the case that $[t,T]$ replaced by $[t,\i)$). Next, for $(t,x)\in[0,T]\times\dbR^n$, we introduce the sets of controls as follows:
\bel{sU^p[t,T]}\ba{ll}
\ns\ds\sU[t,T]=\Bigm\{u:[t,T]\times\O\to\dbR^m\bigm|u(\cd)\hb{ is $\dbF^t$-prograssive measurable}\Big\},\\
\ns\ds\sU^p[t,T]=\Big\{u(\cd)\in\sU[t,T]\bigm|\dbE\(\int_t^T|u(s)|^2ds\)^{p\over2}
<\i\Big\}\equiv L^p_{\dbF^t}(\O;L^2(t,T;\dbR^m)),\qq p\ges1.\ea\ee
For the above state equation, we introduce the following hypothesis:

\ms

{\bf(H1)} It holds $A,C\in\dbR^{n\times n}$ and $B,D\in\dbR^{n\times m}.$

\ms

Under (H1), for any initial pair $(t,x)\in[0,T]\times\dbR^n$ and control $u(\cd)\in\sU^p[t,T]$, there exists a unique solution to the sate equation \rf{state}, denoted by
$$\ba{ll}
\ns\ds X(\cd)=X(\cd\,;t,x,u(\cd))\in L^p_{\dbF^t}(\O;C([t,T];\dbR^n))\equiv\cX^p[t,T]\\
\ns\ds\equiv\Big\{X:[t,T]\times\O\to\dbR^n\bigm|X(\cd)\hb{ is $\dbF^t$-progressively measurable and
has continuous paths, }\\
\ns\ds\qq\qq\qq\qq\qq\qq\qq\dbE\(\sup_{s\in[t,T]}|X(s)|^p\)<\i\Big\}.\ea$$
Moreover, it is standard that the following holds:
\bel{|X|}\dbE\(\sup_{s\in[t,T]}|X(s)|^p\)\les K\[|x|^p+\dbE\(\int_t^T|u(s)|^2ds\)^{p\over2}\].\ee
Further, we have the following stability estimates: For any $t'\in[t,T]$,
\bel{|X-x|}\ba{ll}
\ns\ds\dbE\(\sup_{s\in[t,t']}|X(s)-x|^p\)\les K\dbE\Big\{\[|x|^p+\(\int_t^T|u(s)|^2ds\)^{p\over2}\](t'-t)^{p\over2}+\(\int_t^{t'}
|u(s)|^2ds\)^{p\over2}\Big\},\ea\ee
and if $X'(\cd)=X(\cd\,;t'x',u'(\cd))$ with $(t',x')\in[0,T]\times\dbR^n$, $u'(\cd)\in\sU^p[t',T]$, $t\les t'$. Then
\bel{|X-X|}\ba{ll}
\ns\ds\dbE\(\sup_{s\in[t',T]}|X(s)-X'(s)|^p\)\les K\dbE\[|x'-X(t')|^p+\(\int_t^T|u(s)-u'(s)|^2ds\)^{p\over2}\].\ea\ee
Here, $K$ in \rf{|X|}--\rf{|X-X|} is depending on $p$ and $T$. In what follows, for convenience, we will let $K$ be some generic constant which can be different from line to line. Next, we introduce the following {\it recursive cost functional}:
\bel{J=Y0}J_T(t,x;u(\cd))=Y_T(t),\ee
where $(Y_T(\cd),Z_T(\cd))$ is the adapted solution to following linear backward stochastic differential equation (BSDE, for short):
\bel{BSDE0}\left\{\2n\ba{ll}
\ds dY_T(s)=\(EY_T(s)+FZ_T(s)-f(s,X(s),u(s))\)ds+Z_T(s)dW(s),\qq s\in[t,T],\\
\ns\ds Y_T(T)=\xi_T,\ea\right.\ee
for some non-negative constants $E$ and $F$, with
\bel{f}\left\{\2n\ba{ll}
\ns\ds f(s,x,u)=\lan Qx,x\ran+2\lan Sx,u\ran+\lan Ru,u\ran+2\lan q(s),x\ran+2\lan r(s),u\ran,\\
\ns\ds\qq\qq\qq\qq\qq\qq\qq\qq\qq\qq\qq\qq\qq\qq\qq\qq s\in[t,T],\\
\ns\ds\xi_T=\lan GX(T),X(T)\ran+2\lan g,X(T)\ran,\ea\right.\ee
for some matrices $Q,G$, $S$ and $R$, vector $g$, and deterministic maps $q(\cd)$ and $r(\cd)$. We assume that $q(\cd)$ and $r(\cd)$ to be deterministic, since they are required to be $\dbF^t$-progressively measurable, which implies that $q(t)$ and $r(t)$ are $\dbE_t^t=\cN$ measurable. Thus, $q(t)$ and $r(t)$ are almost sure deterministic, for all $t\in[0,T]$. Hence, it is convenient to assume that they are deterministic. It is clear that by \rf{f}, $s\mapsto f(s,X(s),u(s))\equiv f(s)$ is defined on any finite horizon $[t,T]$, whereas $\xi_T$ is defined only at $T$ (unless $G=0$ and $g=0$). Note that to study the optimal control problem in a finite horizon with recursive cost functional, one has to face the well-posedness of BSDE \rf{BSDE0} (in $L^1$). We see that in general, for $u(\cd)\in\sU^2[t,T]$, the map $f(\cd)\in L^1_{\dbF^t}(\O;L^1(t,T;\dbR))$, and the random variable $\xi_T\in L^1_{\cF_T^t}(\O;\dbR)$. Thus, \rf{BSDE0} is a linear BSDE in the finite horizon with the nonhomogeneous term and the terminal state in $L^1$. It is known that such a BSDE is not well-posed (see \cite{Briand-Hu 2005, Fan-Liu 2010, Fan 2016, Klimsaik-Rzymowski 2024}, among many other papers). Thus, it is a little difficult to discuss such an optimal control problem on $\sU^2[t,T]$.

\ms

The main feature of recursive cost \rf{J=Y0} is that the current cost is depending on the future costs. Such a functional was firstly introduced by Duffie--Epstein in 1992 (\cite{Duffie-Epstein 1992a,Duffie-Epstein 1992b}) (which was called differential utility). Shortly after, people realized that this was merely some special cases of BSDEs. Based on that, in 1992, Peng investigated general nonlinear optimal control problems with recursive cost functional the first time (\cite{Peng 1992}), and introduced the so-called $g$-expectation notion (a generalized expectation, which could be nonlinear, but is time-consistent) in 1997 (\cite{Peng 1997}). The notion of differential utility was extended by Quenez--Lazrak \cite{Lazrak-Quenez 2003}, and Lazrak \cite{Lazrak 2004} in the beginning of 2000s to a general BSDE framework. In 2008, Yong studied LQ problems under (nonlinear) $g$-expectation introduced in \cite{Peng 1997} with the framework of $L^p$ space ($p>1$) \cite{Yong 2008}. Thus, the linear structure was partially ruined (since the expectation was nonlinear) and moreover, it could not cover the case of $p=2$. Because of that, the Riccati type differential equation was derived only for a very special case in \cite{Yong 2008}. Recently, in \cite{Li-Yong 2026}, a general LQ stochastic optimal control with recursive cost in a finite horizon was investigated. Among other things, it was established the open-loop and closed-loop solvability of the optimal control problem and their characterizations in terms of the solvability of a forward-backward stochastic differential equations (FBSDEs, for short) and differential Riccati equation (DRE, for short).

\ms

Now, we may similarly define $\sU[t,\i)$ and $\sU^p[t,\i)$ (compare with \rf{sU^p[t,T]}; we have an infinite horizon here). Then, under (H1), for any initial pair $(t,x)\in[0,\i)\times\dbR^n$ and control $u(\cd)\in\sU^2[t,\i)$, the following controlled system (in infinite horizon):
\bel{state2}\left\{\2n\ba{ll}
\ds dX(s)=[AX(s)+Bu(s)]ds+[CX(s)+Du(s)]dW(s),\q s\in[t,\i),\\
\ns\ds X(t)=x,\ea\right.\ee
admits a unique solution $X(\cd)\equiv X(\cd\,;t,x,u(\cd))$ (defined on $[t,\i)$). Next, we define
\bel{J=Y2}J(t,x;u(\cd))=\widehat Y(t),\ee
where $(\widehat Y(\cd),\widehat Z(\cd))$ is the adapted solution to the following BSDE:
(compare with \rf{BSDE0}; we have an infinite horizon here):
\bel{BSDE1}d\widehat Y(s)=[E\widehat Y(s)+F\widehat Z(s)-f(s,X(s),u(s))]ds+\widehat Z(s)dW(s),\qq s\in[t,\i),\ee
for some non-negative constants $E$ and $F$, where $f(\cd)\equiv f(\cd\,,X(\cd),u(\cd))\in L^1_{\dbF^t}(\O;L^1(t,\i;\dbR))$ is defined as the first in \rf{f} (with $s\in[t,T]$ replaced by $s\in[t,\i)$).
If BSDE \rf{BSDE1} admits a unique adapted solution $(\widehat Y(\cd),\widehat Z(\cd))$, then we can define the corresponding recursive cost functional to be \rf{J=Y2}. Then, we could pose the optimal control problem and study the so-called open-loop and closed-loop solvability and their characterizations in terms of solvability of an FBSDE and that of algebraic Riccati equation (ARE, for short). etc.

\ms

We have seen that even for finite horizon case, the well-posedness of BSDE in $L^1$ is subtle (see \cite{Li-Yong 2026}). In the current paper, we need to consider the BSDE in $L^1$ and in an infinite horizon. In \cite{Peng-Shi 2000,Yin 2008,Sun-Yong 2018,Sun-Yong 2020}, it was investigated the unique solvability in $L^2$ sense. In $L^p$ case, with $p\ges2$, see \cite{Fuhrman-Tessitore 2004}. More recently, the solvability of FBSDE with jumps on infinite horizon and its application to backward linear-quadratic (LQ) problems in \cite{Yu 2017}, the stochastic linear-quadratic optimal control problem with random time-varying coefficients on infinite horizon in \cite{Wei-Yu 2018}, see \cite{Luo-Li-Wei 2025} also for infinite time horizon stochastic recursive control problem with jumps. However, they still require nonhomogeneous term $f(\cd)$ to be at least $L^p, p\ges2$ integrable. To the best of our knowledge, we do not find works to address the solvability of this type of BSDE on infinite horizon when $f(\cd)$ is only $L^1$ integrable.

\ms

In the current paper, to correctly define the recursive cost functional, we need to first look at BSDEs in $L^1$ space with an infinite horizon. With the helps from the results in finite horizons, we obtain the solvability of such BSDEs with the nonhomogeneous term in $L^1$, extending the results of \cite{Sun-Yong 2020}. In our case, due to the fact that the non-homogeneous process is of quadratic form in the state-control pairs, which leads to the weighted $L^2$-stabilizability of the controlled system. Through the characterization, we find weighted $L^2$-stabilizability of the given controlled system is equivalent to the (usual) $L^2$-stabilizability of a modified controlled system. Thus, the assumption of the weighted $L^2$-stabilizability of the control system is reasonable and easy to check. Further, we find that with a proper transformation, the stochastic LQ optimal control problem with recursive cost functional is equivalent to a modified classical LQ problem. Therefore, all the results for the classical LQ problem can be translated. By translation, we obtain the open-loop and closed-loop solvability of the LQ problem with recursive cost functional in an infinite horizon, and their characterizations in terms of the solvability of an FBSDE, and that of ARE, respectively. Finally, we investigate some extensions.

\ms

\section{BSDEs of $L^1$ in an Infinite Horizon}

In this section, we are going to look at the well-posedness of BSDEs \rf{BSDE1} in $L^1$ and in an infinite horizon with $f(\cd\,,X(\cd),u(\cd))$ replaced by a general $f(\cd)$. First of all, for $0\les t<T$, $p,p_1,p_2\ges1$, and $k\ges1$, let
$$\ba{ll}
\ns\ds L^p_{\cF_T^t}(\O;\dbR^k)=\Big\{\xi:\O\to\dbR^k\bigm|\xi\hb{ is $\cF_T^t$-measurable, }\dbE|\xi|^p<\i\Big\},\\
%
\ns\ds L^{p_1}_{\dbF^t}(\O;L^{p_2}(t,T;\dbR^k))=\Big\{f:[t,T]\times\O\to\dbR^k\bigm|f(\cd)
\hb{ is $\dbF^t$-prgressively measurable,}\\
\ns\ds\qq\qq\qq\qq\qq\qq\qq\qq\qq\qq\dbE\(\int_t^T|f(s)|^{p_2}ds\)^{p_1\over p_2}<\i\Big\},\ea$$
The spaces $L^\i_{\cF_T^t}(\O;\dbR^k)$ and $L^\i_{\dbF^t}(\O;L^\i(t,T;\dbR^k))$ can be defined similarly. We now present the following result for the BSDEs \rf{BSDE0} (in finite
horizon and in $L^1$) (see \cite{Li-Yong 2026}, for a more general case with constants $E$ and $F$ replaced by some maps).

\bl{dm} \sl Let $\m(\cd)$ solve the following one-dimensional SDE:
\bel{dF}\left\{\2n\ba{ll}
\ds d\m(s)=-E\m(s)ds-F\m(s)dW(s),\qq s\ges0,\\
\ns\ds\m(0)=1.\ea\right.\ee
Then
\bel{m}\m(s)=e^{-Es-{1\over2}F^2s-FW(s)},\qq s\ges t.\ee
Consequently,
\bel{Em}\dbE[\m(s)]=e^{-Es},\qq s\ges0.\ee
Moreover, for any $\xi_T\in L^\i_{\cF_T^t}(\O;\dbR^n)$, and $f(\cd)\in L^\i_{\dbF^t}(\O;L^\i(t,T;\dbR))$, BSDE \rf{BSDE0} admits a unique adapted solution $(Y_T(\cd),Z_T(\cd))$ and the following representation holds:
\bel{Y_T}Y_T(s)=\dbE_s^t\[{\m(T)\over\m(s)}\xi_T+\int_s^T{\m(\t)\over\m(s)}f(\t)d\t\],\qq s\in[t,T].\ee
Conversely, if \rf{Y_T} holds for some
\bel{f-p}\xi_T\in L^\i_{\cF_T^t}(\O;\dbR^n),\qq f(\cd)\in L^\i_{\dbF^t}(\O;L^\i(t,T;\dbR)),\ee
then for some $Z_T(\cd)\in L^p_{\dbF^t}(\O;L^2(t,T;\dbR^n))$, (with $p>1$) $(Y_T(\cd),Z_T(\cd))$ is the adapted solution to BSDE \rf{BSDE0}.

\el

\it Proof. \rm Note that SDE \rf{dF} is a classical SDE, it is well-posed. Thus we may let $\m(\cd)$ be the solution of \rf{dF}. Consequently, by It\^o's formula,
$$d\(\ln\m(s)\)={1\over\m(s)}\[E\m(s)ds+F\m(s)dW(s)\]-{1\over2}
{1\over\m(s)^2}F^2\m(s)^2ds=\(E-{1\over2}F^2\)ds+FdW(s).$$
Hence, \rf{m} follows. Now, for some $\xi_T\in L^\i_{\cF_T^t}(\O;\dbR^k)$ and $f(\cd)\in L^\i_{\dbF^t}(\O;L^\i(t,T;\dbR))$, satisfying \rf{f-p}, BSDE \rf{BSDE0} admits a unique adapted solution $(Y_T(\cd),Z_T(\cd))$, and
$$\dbE\[\(\sup_{s\in[t,T]}|Y_T(s)|^p\)+\(\int_t^T|Z_T(s)|^2ds\)^{p\over2}\]\les K\dbE\[|\xi_T|^p+\(\int_t^T|f(s)|ds\)^p\],$$
for some $p>1$ and $K>0$. By It\^o's formula, we have
$$\ba{ll}
\ns\ds d\(\m(s)Y_T(s)\)=\(-E\m(s)ds-F\m(s)dW(s)\)Y_T(s)\\
\ns\ds\qq\qq\qq\qq+\m(s)
\[\(EY_T(s)+FZ_T(s)-f(s)\)ds+Z_T(s)dW(s)\]+F\m(s,t)Z_T(s)ds\\
\ns\ds\qq\qq\qq\q=-\m(s)f(s)ds+\m(s)\(-FY_T(s)+Z_T(s)\)dW(s),\qq s\in[t,T].\ea$$
Thus, we have
$$\ba{ll}
\ns\ds\dbE_s^t\[\m(T)\xi_T-\m(s)Y_T(s)\]=-\dbE_s^t\int_s^T
\m(\t)f(\t)d\t\ea$$
The rest conclusions follow easily. \endpf

\ms

Note that
$${\m(s)\over\m(t)}=e^{-E(s-t)-{F^2\over2}(s-t)-F[W(s)-W(t)]},\qq s\in[t,T].$$
Thus, ${\m(\cd)\over\m(t)}$ is $\dbF^t$-progressively measurable. It is suggested by \rf{Y_T} that we should introduce the following: for any $0\les t<T$, $p,p_1,p_2\ges1$ and $k\ges1$,
$$\ba{ll}
\ns\ds L^{\m,p}_{\cF_T^t}(\O;\dbR^k)\equiv\Big\{\xi:\O\to\dbR\bigm|\xi\hb{ is $\cF_T^t$-measurable, }\dbE\[{\m(\t)\over\m(t)}|\xi|^p\]<\i\Big\},\\
\ns\ds L^{p_1}_{\dbF^t}(\O;L^{\m,p_2}(t,T;\dbR^k))\equiv\Big\{f:\O\times[t,T]\to\dbR^k\bigm|
f(\cd)\hb{ is $\dbF^t$-prograssively measurable,}\\
\ns\ds\qq\qq\qq\qq\qq\qq\qq\qq\qq\qq\qq\dbE\(\int_t^T{\m(s)\over\m(t)}|f(s)|^{p_2}ds\)^{p_1\over p_2}<\i\Big\}.\ea$$
Any element in the above spaces is said to be weighted integrable. Next, we write
\bel{Y_T*}Y_T(s)=\dbE_s^t\[{\m(T)\over\m(s)}\(\xi_T^+-\xi^-_T\)+\int_s^T{\m(\t)\over
\m(s)}\(f(\t)^+-f(\t)^-\)d\t\],\qq s\in[t,T],\ee
where $a^+=\max\{a,0\}$ and $a^-=\max\{-a,0\}$ for any $a\in\dbR$. Then the above is well-defined (including the definite $\pm\i$) when
$$\ba{ll}
\ns\ds\xi^-_T\in L^{\m,1}_{\cF_T^t}(\O;\dbR),\qq f(\cd)^-\in L^1_{\dbF^t}(\O;L^{\m,1}(t,T;\dbR)),\ea$$
which is the case that $Y_T(s)\les\i$, for some $s\in[t,T]$ (the indefinite form $\i-\i$ is avoided); and
$$\xi^+_T\in L^{\m,1}_{\cF_T^t}(\O;\dbR),\q f(\cd)^+\in L^1_{\dbF^t}(\O;L^{\m,1}(t,T;\dbR)),$$
which is the case that $Y_T(s)\ges-\i$, for some $s\in[t,T]$ (the indefinite form $(-\i)-(-\i)$ is avoided). This can be stated and proved (with small modification) for $E$ and $F$ replaced by functions (instead of constants).  Hence, in some sense, we extend the results of \cite{Peng 1997}.

\ms

Now, we take $G=0$ and $g=0$ (so that $\xi_T=0$). Clearly, in this case,
\bel{Y_T**}Y_T(s)=\dbE_s^t\int_s^T{\m(\t)\over\m(s)}
\(f(\t)^+-f(\t)^-\)d\t,\qq s\in[t,T].\ee
To study the above when $T\to\i$, it is natural to define $L^1_{\dbF^t}(\O;L^{\m,1}(t,\i;\dbR))$ (replacing $[t,T]$ by $[t,\i)$).
The limit in \rf{Y_T**} will exist, including the definite forms of $\pm\i$, provided
\bel{mf+}f(\cd)^+\hb{ or }f(\cd)^- \in L^1_{\dbF^t}(\O;L^{\m,1}(t,\i;\dbR)).\ee
We now denote
$$\h L^1_{\dbF^t}(\O;L^{\m,1}(t,\i;\dbR))=\Big\{f:\O\times[t,\i)\to\dbR\bigm|f(\cd)^+\hb{ or }f(\cd)^-\in L^1_{\dbF^t}(\O;L^{\m,1}(t,\i;\dbR))\Big\}.$$
Clearly, $L^1_{\dbF^t}(\O;L^{\m,1}(t,\i;\dbR))$ is dense in $\h L^1_{\dbF^t}(\O;L^{\m,1}(t,\i;\dbR))$. But the latter is not a subspace, since
$$f(\cd)=f(\cd)^+-f(\cd)^-\in L^1_{\dbF^t}(\O;L^{\m,1}(t,\i;\dbR))\q\Ra\q f(\cd)^+\q\&\q f(\cd)^-\in L^1_{\dbF^t}(\O;L^{\m,1}(t,\i;\dbR)),$$
but
$$f(\cd)=f(\cd)^+-f(\cd)^-\in\h L^1_{\dbF^t}(\O;L^{\m,1}(t,\i;\dbR))\q\Ra\q f(\cd)^+\q\hb{ or }\q f(\cd)^-\in L^1_{\dbF^t}(\O;L^{\m,1}(t,\i;\dbR)).$$

\ms

From \rf{Y_T}, we see that for $f(\cd)\in\h L^1_{\dbF^t}(\O;L^{\m,1}(0,\i;\dbR))$, then \rf{Y_T**} is well-defined (it might be the infinity $\pm\i$). In particular, the recursive cost functional \rf{J=Y2} could be well-defined if $f(\cd)\in L^1_{\dbF^t}(\O;L^{\m,1}(0,\i;\dbR))$. For the weight function $\m(\cd)$, we know that ${\m(\cd)\over\m(t)}\in L^1_{\dbF^t}(\O;L^1(t,\i;\dbR))$ and
\bel{Em}\dbE\int_t^\i{\m(\t)\over\m(t)}d\t=\int_t^\i e^{-E(s-t)}ds={1\over E},\ee
provided $E>0$. But ${\m(\cd)\over\m(t)}\notin L^1_{\dbF^t}(\O;L^1(t,\i;\dbR^n))$, if $E\les 0$. Next, for $f(\cd)\in L^{p_1}_{\dbF^t}(\O;L^{p_2}(t,\i;\dbR))$,
\bel{Em*}\ba{ll}
\ns\ds\dbE\Big|\int_t^\i{\m(\t)\over\m(t)}f(\t)d\t\Big|\les\dbE\Big\{\[\int_t^\i
\({\m(\t)\over\m(t)}\)^{p_2\over p_2-1}d\t\]^{p_2-1\over p_2}\(\int_t^\i|f(\t)|^{p_2}d\t\)^{1\over p_2}\Big\},\ea\ee
provided,
$${(p_2-1)p_1\over p_2(p_1-1)}\ges1\q\iff\q p_2\ges p_1.$$
In this case,
$$\ba{ll}
\ns\ds\int_t^\i\dbE\({\m(\t)\over\m(t)}\)^{p_2\over p_2-1}d\t=\int_t^\i \dbE\[e^{(-E(\t-t)-{F^2\over2}(\t-t)-F(W(\t)-W(t))){p_2\over p_2-1}}\]d\t\\
\ns\ds=\int_t^\i\[e^{-E{p_2\over{p_2-1}}(\t-t)-{F^2\over2}{p_2\over {p_2-1}}(\t-t)+{1\over2}({Fp_2\over{p_2-1}})^2(\t-t)}\]d\t=\int_t^\i e^{-{p_2\over{p_2-1}}(E-{1\over2}{1\over p_2-1}F^2)(\t-t)}d\t.\ea$$
Hence, for $f(\cd)\in L^{p_1}_{\dbF^t}(\O;L^{p_2}(t,\i;\dbR))$ with $p_2\ges p_1$, $\int_t^\i\({\m(\t)\over\m(t)}\)f(\t)d\t\in L^1_{\cF^t_\i}(\O;\dbR)$ if $E>0$ and
\bel{p_2}E-{1\over2(p_2-1)}F^2>0\q\iff\q p_2>1+{F^2\over2E}.\ee
Recall that for BSDE \rf{BSDE1}, it is known from \cite{Sun-Yong 2020} that for $f(\cd)\in L^2_{\dbF^t}(\O;L^2(t,\i;\dbR))$, it admits a unique adapted solution $(Y(\cd),Z(\cd))$ if the system $[-E,-F]$ is stable, or, equivalently,
\bel{EF}E-{F^2\over2}>0\q\iff\q 2>1+{F^2\over2E}.\ee
Clearly, \rf{p_2} is a generalization of \rf{EF}. With a little modifications, one can claim that for any $f(\cd)\in L^{p_1}_\dbF(\O;L^{p_2}(0,\i;\dbR^n))$, with \rf{p_2}, BSDE \rf{BSDE1} has a unique adapted solution $(Y(\cd),Z(\cd))$, extending that in \cite{Sun-Yong 2020}.

\section{Weighted $L^2$-Stabilizability and Optimal Control Problems with Recursive Cost Functionals}

In what follows, we let $\dbS^k$, $\dbS_+^k$, $\dbS_{++}^k$ be the sets of all $(k\times k)$ symmetric matrices, positive semi-definite matrices and positive definite matrices, respectively.
In order the recursive cost functional to be well-defined, we observe the following:
$$\ba{ll}
\ns\ds|\widehat Y(s)|=\Big|\dbE_s^t\int_s^\i{\m(\t)\over\m(s)}\[\lan\begin{pmatrix}Q&S^\top\\ S&R\end{pmatrix}\begin{pmatrix}X(\t)\\ u(\t)\end{pmatrix},\begin{pmatrix}X(\t)\\ u(\t)\end{pmatrix}\ran+2\lan \begin{pmatrix}q(\t)\\ r(\t)\end{pmatrix},\begin{pmatrix}X(\t)\\ u(\t)\end{pmatrix}\ran\]d\t\Big|\\
\ns\ds\les\dbE_s^t\[\int_s^\i{1\over\m(s)}\Big|\lan\begin{pmatrix}Q&S^\top\\ S&R\end{pmatrix}\begin{pmatrix}\sqrt{\m(\t)}X(\t)\\ \sqrt{\m(\t)}u(\t)\end{pmatrix},\begin{pmatrix}\sqrt{\m(\t)}X(\t)\\ \sqrt{\m(\t)}u(\t)\end{pmatrix}\ran+2\lan\begin{pmatrix}\sqrt{\m(\t)}q(\t)\\ \sqrt{\m(\t)}r(\t)\end{pmatrix},\begin{pmatrix}\sqrt{\m(\t)}X(\t)\\ \sqrt{\m(\t)}u(\t)\end{pmatrix}\ran\Big|d\t\]\\
\ns\ds\les\dbE_s^t\int_s^\i\[\(|Q|+2|S|+|R|+1\){\m(\t)\over\m(s)}\(|X(\t)|^2+|u(\t)|^2\)
+{\m(\t)\over\m(s)}\(|q(\t)|^2
+|r(\t)|^2\)\]d\t.\ea$$
Hence,
\bel{|J|}\ba{ll}
\ns\ds\Big|J(t,x;u(\cd))\Big|\les \(|Q|+2|S|+|R|+1\)\dbE\int_t^\i{\m(\t)\over\m(t)}\(|X(\t)|^2+|u(\t)|^2\)d\t\\
\ns\ds\qq\qq\qq\qq\qq+\dbE
\int_t^\i{\m(\t)\over\m(t)}\(|q(\t)|^2+|r(\t)|^2\)d\t,\ea\ee
provided the right-hand side is finite. Let $k\ges 1$, and
$$\ba{ll}
\ns\ds L^{\m,2}(0,\i;\dbR^k)\equiv \left\{\varphi: [0,\i)\to\dbR^k \Big|\varphi\,\, \text{is}\,\, \text{Lebesgue measurable and}\,\int_0^\i \m(t)|\varphi(t)|^2dt < \infty\right\}. \\
\ea $$

In order for the second integral term on the right-hand side of the above to be defined, we had better to introduce the following hypothesis.

\ms

{\bf(H2)} Let $E,F>0$, $Q\in\dbS^n$, $S\in\dbR^{m\times n}$, $R\in\dbS^n$, $q(\cd)\in L^{\m,2}(0,\i;\dbR^n)$ and $r(\cd)\in L^{\m,2}(0,\i$; $\dbR^m)$.

\ms

In order for the first integral term of the right-hand side of \rf{|J|} makes sense, we need $X(\cd)\in L^2_{\dbF^t}(\O;$ $L^{\m,2}(t,\i;
\dbR^n))$ under control $u(\cd)\in\sU^{\m,2}[t,\i)\equiv L^2_{\dbF^t}(\O;L^{\m,2}(t,\i;\dbR^m))$. To approach this, we take the following method. With control $u(\cd)=0$, the uncontrolled system is denoted by $[A,C]$. Hence, for initial pair $(t,x)\in[0,\i)\times\dbR^n$ and $u(\cd)=0\in\sU^{\m,2}[t,\i)$, the solution $X(\cd\,;t,x,0)$ to the uncontrolled system $[A,C]$ looks like the following:
\bel{state3}\left\{\2n\ba{ll}
\ds dX(s;t,x,0)=AX(s;t,x,0)ds+CX(s;t,x,0)dW(s),\q s\in[t,\i),\\
\ns\ds X(t;t,x,0)=x,\ea\right.\ee
Now, let us introduce the following definition.

\bde{} (i) System $[A,C]$ is said to be weighted $L^2$-stable if for any initial pair $(t,x)\in[0,\i)\times\dbR^n$, the solution $X(\cd\,;t,x,0)$ of system $[A,C]$ is in $L^2_\dbF(\O;L^{\m,2}(t,\i;\dbR))$, i.e.,
$$\dbE\int_t^\i\m(\t)|X(\t;t,x,0)|^2d\t<\i.$$

(ii) System $[A,C;B,D]$ is said to be weighted $L^2$-stabilizable if there exists a $\Th\in\dbR^{n\times m}$ such that $[A+B\Th,C+D\Th]$ is weighted
$L^2$-stable. Any such a $\Th$ is called a weighted $L^2$-stabilizer of system $[A,C;B,D]$.

\ede

In what follows, we let
\bel{cX^{m,2}}\cX^{\m,2}[t,\i)=\[\bigcup_{T>t}L^2_{\dbF^t}\O;C([t,T];\dbR^n))\]
\bigcap L^2_{\dbF^t}(\O;L^{\m,2}(t,\i;\dbR^n)),\ee
and
\bel{cX^2}\cX^2[t,\i)=\[\bigcup_{T>t}L^2_{\dbF^t}(\O;C([t,T];\dbR^n))\]\bigcap L^2_{\dbF^t}(\O;L^2(t,\i;\dbR^n)).\ee

Clearly, $[A,C]$ is weighted $L^2$-stable if and only if
$$X(\cd\,;t,x,0)\in\cX^{\m,2}[t,\i),\qq(t,x)\in[0,\i)\times\dbR^n,$$
We denote $\hTh^\m[A,C;B,D]$ to be the set of all weighted $L^2$-stabilizers of $[A,C;B,D]$. Now, for any {\it state feedback strategy } $(\Th,v(\cd))\in\hTh^\m[A,C;B,D]\times\sU^{\m,2}[t,\i)$, we may consider the following system:
\bel{closed-loop}\left\{\2n\ba{ll}
\ds dX(s)=\[A^\Th X(s)+Bv(s)\]ds+\[C^\Th X(s)+Dv(s)\]dW(s),\\
\ns\ds X(t)=x,\ea\right.\ee
where
\bel{A^Th}A^\Th=A+B\Th,\qq C^\Th=C+D\Th.\ee
We call \rf{closed-loop} the {\it closed-loop system} under strategy $(\Th,v(\cd))$. The solution of \rf{closed-loop} is denoted by $X(\cd)\equiv X(\cd\,;t,x,(\Th,v(\cd)))$ $\in\cX^{\m,2}[t,\i)$, which  depends on the initial pair $(t,x)\in[0,\i)\times\dbR^n$. It is easy for us to realize that \rf{closed-loop} is equivalent to the controlled system \rf{state} under the control
\bel{control}u(s)=\Th X(s)+v(s),\qq s\in[t,\i),\ee
with $X(\cd)$ being the corresponding state process. Thus,
\bel{X=X}X(\cd\,;t,x,u(\cd))=X(\cd\,;t,x,(\Th,v(\cd))).\ee
Because of \rf{control} (and \rf{X=X}), we call such a control an {\it outcome} of the strategy $(\Th,v(\cd))$, which depends on the initial pair $(t,x)$ (since both sides of \rf{X=X} do).

\ms

We now try to get the characterization of weighted $L^2$-stability of uncontrolled system $[A,C]$ and weighted $L^2$-stabilizability of controlled system $[A,C;B,D]$. To this end, let $X_0(\cd)=X(\cd\,;t,x,0)$. Then, for any $P\in\dbS^n$ ($P=I$ is allowed), we have (by It\^o's formula)
$$\ba{ll}
\ns\ds d\lan PX_0(s),X_0(s)\ran=\(\lan PAX_0(s)ds+PCX_0(s)dW(s),X_0(s)\ran+\lan PX_0(s),AX_0(s)ds+CX_0(s)dW(s)\ran\\
\ns\ds\qq\qq\qq\qq\qq+\lan PCX_0(s),CX_0(s)\ran ds\\
\ns\ds\qq\qq\qq\q=\lan(PA+A^\top P+C^\top PC)X_0(s),X_0(s)\ran ds+\lan(PC+C^\top P)X_0(s),X_0(s)\ran dW(s).\ea$$
Thus
$$\ba{ll}
\ns\ds d\(\m(s)\lan PX_0(s),X_0(s)\ran\)=\(-E\m(s)ds-F\m(s)dW(s)\)\lan PX_0(s),X_0(s)\ran\\
\ns\ds\qq\qq\qq\qq\qq\qq+\m(s)\(\lan(PA+A^\top P+C^\top PC)X_0(s),X_0(s)\ran ds\\
\ns\ds\qq\qq\qq\qq\qq\qq+\lan(PC+C^\top P)X_0(s),X_0(s)\ran dW(s)\)-F\m(s)\lan(PC+C^\top P)X_0(s),X_0(s)\ran ds\\
\ns\ds\qq\qq\qq\qq\qq\q=\lan\m(s)[-EP-F(PC+C^\top P)+PA+A^\top P+C^\top PC]X_0(s),X_0(s)\ran ds\\
\ns\ds\qq\qq\qq\qq\qq\qq+\m(s)\lan(-FP+PC+C^\top P)X_0(s),X_0(s)\ran dW(s).\ea$$
Consequently,
\bel{m|X|^2}\ba{ll}
\ns\ds\dbE\[\m(s)\lan PX_0(s),X_0(s)\ran-\m(t)\lan Px,x\ran\]\\
\ns\ds=\dbE\int_t^s\m(\t)\lan[-EP
-F(PC+C^\top P)+PA+A^\top P+C^\top PC]X_0(\t),X_0(\t)\ran d\t.\ea\ee
From the above, we have the following result.

\bp{P*}\label{prop3.2} \sl System $[A,C]$ is weighted $L^2$-stable if and only if there exists a $P\in\dbS_{++}^n$ such that the following Lyapunov inequality holds
\bel{Lyapunov1}-EP-F(PC+C^\top P)+PA+A^\top P+C^\top PC<0.\ee
Consequently, for any $\L\in\dbS^n_{++}$, the Lyapunov equation
\bel{Lyapunov2}-EP-F(PC+C^\top P)+PA+A^\top P+C^\top PC+\L=0.\ee
admits a unique solution $P\in\dbS^n_{++}$ given by
\bel{P}P=\dbE\int_0^\i\m(s)\Psi(s,0)^\top\L\Psi(s,0)ds,\ee
where
$$\left\{\2n\ba{ll}
\ns\ds d\Psi(s,t)=A\Psi(s,t)ds+C\Psi(s,t)dW(s),\qq s\ges t,\\
\ns\ds\Psi(t,t)=I.\ea\right.$$

\ep

\it Proof. \rm Necessity. For any fixed $\L\in\dbS^n$, consider the linear ODE on $[0,\i):$
\bel{P}\left\{\2n\ba{ll}
\ns\ds\dot\Th(s)=-E\Th(s)-F[\Th(s)C+C^\top \Th(s)]+\Th(s)A+A^\top\Th(s)+C^\top\Th(s)C+\L=0, \\
\ns\ds\Th(0)=0.\ea\right.\ee
Clearly, it has a unique solution $\Th(s)$ defined on $[0,\i)$. For any fixed $\t>0$, the function
$$\Th^\t(s):=\Th(\t-s),\quad s\in[0,\t]$$
solves the equation
\bel{2.14}\left\{\2n\ba{ll}
\ns\ds\dot{\Th}^\t(s)-\1n E\Th^\t(s)\1n-\1n F[\Th^\t(s)C\1n+\1n C^\top\Th^\t(s)]\1n+\1n\Th^\t(s)A\1n+\1n A^\top\Th^\t(s)\1n+\1n C^\top\Th^\t(s)C\1n+\1n\L\1n=
\1n0,\q s\in[0,\t],\\
\ns\ds\Th^\t(\t)=0.\ea\right.\ee
Let $X_0(\cd)\equiv X(\cd\,;0,x)$ be the solution to system $[A,C]$ with initial pair $(0,x)$ and note that $X_0(s)$ admits the representation $X_0(s)=\Psi(s,0)x$. Thus, we have
$$\ba{ll}
\ns\ds d\(\Th^\t(s)X_0(s)\)=\[\dot\Th^\t(s)X_0(s)+\Th^\t(s)AX_0(s)\]ds
+\Th^\t(s)CX_0(s)dW(s)\ea$$
Then
$$\ba{ll}
\ns\ds d\(\lan\Th^\t(s)X_0(s),X_0(s)\ran\)=\lan[\dot\Th^\t(s)X_0(s)+\Th^\t(s)AX_0(s)]ds
+\Th^\t(s)CX_0(s)dW(s),X_0(s)\ran\\
\ns\ds\qq\qq\qq\qq\qq\qq+\lan\Th^\t(s)X_0(s),AX_0(s)ds+CX_0(s)dW(s)\ran+\lan\Th^\t(s)CX_0(s),CX_0(s)\ran ds\\
\ns\ds\qq\qq\qq\qq\qq\q=\lan[\dot\Th^\t(s)+\Th^\t(s)A+A^\top
\Th^\t(s)+C^\top\Th^\t(s)C]X_0(s),X_0(s)\ran ds\\
\ns\ds\qq\qq\qq\qq\qq\qq+\lan[C^\top\Th^\t(s)+\Th^\t(s)C]X_0(s),X_0(s)\ran dW(s)\ea$$
Hence,
\bel{}\ba{ll}
\ns\ds-\lan\Th^\t(0)x,x\ran=\dbE\[\m(\t)\lan\Th^\t(\t)X(\t),X(\t)\ran-\lan\Th^\t(0) x,x\ran\]=\dbE\int_0^\t d\(\m(s)\lan\Th^\t(s)X_0(s),X_0(s)\ran\)\\
\ns\ds\qq=\dbE\int_0^\t\(-E\m(s)ds-F\m(s)dW(s)\)\lan\Th^\t(s)
X_0((s),X_0(s)\ran\\
\ns\ds\qq\qq\qq\qq+\m(s)\lan[\dot\Th^\t(s)+\Th^\t(s)A+A^\top
\Th^\t(s)+C^\top\Th^\t(s)C]X_0(s),X_0(s)\ran\\
\ns\ds\qq\qq\qq\qq-F\m(s)\lan\Th^\t(s)X_0(s),X_0(s)\ran\)ds\\
\ns\ds\qq=\dbE\int_0^\t\m(s)\lan[\dot{\Th}^\t(s)+\Th^\t(s)A+A^\top\Th^\t(s)
+C^\top\Th^\t(s)C-E\Th^\t(s)-F\Th^\t(s)]X_0(s),X_0(s)\ran\\
\ns\ds\qq=-\dbE\int_0^\t\m(s)\lan\L X_0(s),X_0(s)\ran ds=-x^\top\[\dbE\int_0^\t\m(s)\Psi(s)^\top\L\Psi(s)ds\]x.\ea\ee
It follows that
\bel{}\Th(\t)=\Th^\t(0)=\dbE\int_0^\t\m(s)\Psi(s,0)^\top\L\Psi(s,0)ds,\qq\t\ges0.\ee
If the system $[A, C]$ is weighted $L^2$-stable, one has the following limit:
\bel{}\lim_{\t\to\i}\Th(\t)=\dbE\int_0^\i\m(s)\Psi(s,0)^\top\L\Psi(s,0)ds \equiv P.\ee
Because $\Th(s)$ is the solution to equation \rf{2.14}, we have for any $s>0$,
\bel{}\ba{ll}
\ns\ds\Th(s+1)-\Th(s)=-E\(\int_s^{s+1}\Th(\t)d\t\)
-F\[\(\int_s^{s+1}\Th(\t)d\t\)C+C^\top\(\int_s^{s+1}\Th(\t)d\t\)\]\\
\ns\ds\qq\qq\qq\qq\qq+\(\int_s^{s+1}\Th(\t)d\t\)A+A^\top\(\int_s^{s+1} \Th(\t) d\t\)+C^\top\(\int_s^{s+1}\Th(\t)d\t\)C+\L\ea\ee
Letting $t\to\i$, we obtain \rf{Lyapunov2}.

\ms

Sufficiency. Suppose $P\in\dbS_{++}^n$ satisfies \rf{Lyapunov2}. By Itô's formula, we have for any $s>t$,
$$\ba{ll}
\ns\ds\dbE^t_s\[\lan\m(s)PX_0(s),X_0(s)\ran-\m(t)\lan Px,x\ran\]\\
\ns\ds\qq=\dbE^t_s\int_t^s\m(\t)\lan(-EP-F(PC+C^\top P)+P A+A^\top P+C^\top PC)X_0(\t),X_0(\t)\ran d\t\\
\ns\ds\qq\les-\l\dbE\int_0^s\m(\t)|X_0(\t)|^2d\t,\ea$$
where $-\l<0$ is the largest eigenvalue of $-EP-F(PC+C^\top P)+PA+A^\top P+C^\top PC<0$. Then
$$\l\dbE\int_0^s\m(\t)|X_0(\t)|^2ds\les\lan Px,x\ran-\dbE\lan\m(s)PX_0(s),X_0(s)\ran\les\lan Px,x\ran,\qq\forall s>0.$$
which implies the weighted $L^2$-stability of $[A,C]$. \endpf

\ms

Now, we can write \rf{Lyapunov2} as follows
\bel{Lyapunov3}P\(A-{F\over2}C-{4E+F^2\over8}I\)+\(A-{F\over2}C-{4E+F^2\over8}I\)^\top P+\(C-{F\over2}I\)^\top P\(C-{F\over2}I\)+\L=0.\ee
Thus, $[A,C]$ is weighted $L^2$-stable if and only if $[A-{F\over2}C-{4E+F^2\over8}I,C-{F\over2}I]$ is $L^2$-stable (see \cite{Sun-Yong 2020}).

\ms

Now, having the above Proposition \ref{P*} for uncontrolled system $[A,C]$, it is easy for us to get the following result for controlled system $[A,C;B,D]$.

\bc{} \sl Let the controlled system $[A,C;B,D]$ be given. Then $\Th\in\hTh^\m[A,C;B,D]$ if and only if for each $
\L\in\dbS_{++}^n$, the following admits a solution $P\in\dbS_{++}^n$:
\bel{Lyapunov-Th}-EP-F[P(C+D\Th)+(C+D\Th)^\top P]+P(A+B\Th)+(A+B\Th)^\top P+(C+D\Th)^\top P(C+D\Th)+\L=0.\ee

\ec

%
%
%
%
%

Thus, by writing \rf{Lyapunov-Th} as
\bel{Lyapunov-Th*}\ba{ll}
\ns\ds P\[A-{F\over2}C-{4E+F^2\over8}I+\(B-{F\over2}D\)\Th\]
+\[A-{F\over2}C-{4E+F^2\over8}I+\(B-{F\over2}D\)\Th\]^\top P\\
\ns\ds+\(C-{F\over2}I+D\Th\)^\top P\(C-{F\over2}I+D\Th\)+\L=0,\ea\ee
%
%
%
%
%
%
%
we see that system $[A,C;B,D]$ is weighted $L^2$-stabilizable if and only if $[A-{F\over2}C-{4E+F^2\over8}I,C-{F\over2}I;B-{F\over2}D,D]$ is $L^2$-stabilizable. This means by assuming $\hTh^\m[A,C;B,D]\ne\varnothing$ is reasonably easy to check.

\ms

By the definition, we see that $[A,C]$ is weighted $L^2$-stable if
\bel{stable}-E-F(C+C^\top)+A+A^\top+C^\top C<0,\ee
and $[A,C;B,D]$ is weighted $L^2$-stabilizable if for some $\Th\in\dbR^{m\times n}$ such that
\bel{stabizable}-E-F[(C+D\Th)+(C+D\Th)^\top]+(A+B\Th)+(A+B\Th)^\top+
(C+D\Th)^\top(C+D\Th)<0.\ee

\ms

Now let us look at the nonhomogeneous system
\bel{nonhom}\left\{\2n\ba{ll}
\ds dX(s)=[AX(s)+\varphi(t)]ds+[CX(s)+\rho(t)]dW(s),\q s\in[t,\i),\\
\ns\ds X(t)=x,\ea\right.\ee
With Proposition 3.2, we have the following result for the above nonhomogenous system.
\bp{} \sl Let {\rm(H1)--(H2)} hold. Let $(t,x)\in[0,\i)\times\dbR^n$. Then, for any $\varphi, \rho \in L^2_{\dbF^t}(\O;L^{\m,2}(t,\i;\dbR^n)) $, the solution $X(\cdot) := X(\cdot ; x, \varphi, \rho)$ of (\ref{nonhom}) is in $L^2_{\dbF^t}(\O;L^{\m,2}(t,\i;\dbR^n))$. Moreover, there exists a constant $K>0$, independent of $t, x, \varphi$ and $\rho$, such that
\begin{equation}
    \mathbb{E} \int_t^{\infty}\mu(s)|X(s)|^2 d s \leqslant K\left\{|x|^2+\mathbb{E} \int_t^{\infty}\mu(s)\left[|\varphi(s)|^2+|\rho(s)|^2\right] d s\right\} .
\end{equation}
\ep

\it Proof.
    Since $[A, C]$ is weighted $L^2$-stable, by Proposition 3.2, there exists a $P>0$ such that
\begin{equation}
    P A+A^{\top} P+C^{\top} P C-EP-FPC-FC^\top P+I_n=0 .
\end{equation}
Applying Itô's formula to $s \mapsto\langle \mu(s)P X(s), X(s)\rangle$, we obtain for all $\tau>t$,
\begin{equation}
    \begin{aligned}
& \mathbb{E}\langle \mu(\tau)P X(\tau), X(\tau)\rangle-\langle P x, x\rangle \\
& =\mathbb{E} \int_t^\tau \mu(s)\left[\left\langle\left(P A+A^{\top} P+C^{\top} P C-EP-FPC-FC^\top P\right) X(s), X(s)\right\rangle\right. \\
& \left.\quad+2\left\langle P \varphi(s)+C^{\top} P \rho(s)-FP\rho(s), X(s)\right\rangle+\langle P \rho(s), \rho(s)\rangle\right] d s \\
& =\mathbb{E} \int_t^\tau\mu(s)[-|X(s)|^2+2\left\langle P \varphi(s)+C^{\top} P \rho(s)-FP\rho(s), X(s)\right\rangle +\langle P \rho(s), \rho(s)\rangle] d s
\end{aligned}
\end{equation}
Let $\lambda>0$ be the smallest eigenvalue of $P$ and set
\begin{equation}
    \alpha(s)=P \varphi(s)+C^{\top} P \rho(s)-FP\rho(s), \quad \beta(s)=\langle P \rho(s), \rho(s)\rangle ; \quad s>0 .
\end{equation}
Then by the Cauchy-Schwarz inequality, we have
\begin{equation}
    \begin{aligned}
\lambda \mathbb{E}\mu(\tau)|X(\tau)|^2 & \leqslant \mathbb{E}\langle \mu(\tau)P X(\tau), X(\tau)\rangle \\
& \leqslant\langle P x, x\rangle+\mathbb{E} \int_t^\tau\mu(s)\left[-|X(s)|^2+2\langle\alpha(s), X(s)\rangle+\beta(s)\right] d s \\
& \leqslant\langle P x, x\rangle+\mathbb{E} \int_t^\tau\mu(s)\left[-\frac{1}{2}|X(s)|^2+2|\alpha(s)|^2+\beta(s)\right] d s \\
& =\langle P x, x\rangle+\int_t^\tau\mu(s)\left[-\frac{1}{2} \mathbb{E}|X(s)|^2+2 \mathbb{E}|\alpha(s)|^2+\mathbb{E} \beta(s)\right] d s
\end{aligned}
\end{equation}
It follows from Gronwall's inequality that
\begin{equation}
    \lambda \mathbb{E}\mu(\tau)|X(\tau)|^2 \leqslant\langle P x, x\rangle e^{-(2 \lambda)^{-1} \tau}+\int_t^\tau e^{-(2 \lambda)^{-1}(\tau-s)}\mu(s)\left[2 \mathbb{E}|\alpha(s)|^2+\mathbb{E} \beta(s)\right] d s.
\end{equation}
Finally, integrating $\mathbb{E}\mu(\tau)|X(\tau)|^2$ over $[t, \infty)$, switching the order of integration, we get the desired conclusion.
\endpf

Now, we introduce the following assumption:

\ms

{\bf(H3)} System $[A,C;B,D]$ is weighted $L^2$-stabilizable, i.e., $\hTh^\m[A,C;B,D]\ne\varnothing$.

\ms

Then, the following result holds.

\bp{} \sl Let {\rm(H1)--(H3)} hold. Let $(t,x)\in[0,\i)\times\dbR^n$. Then,
\bel{sU^{m,2}}\ba{ll}
\ns\ds\sU^{\m,2}[t,\i)\subseteq\Big\{\Th X(\cd)+v(\cd)\bigm|(\Th,v(\cd))\in\hTh^\m[A,C;B,D]\times\sU^{\m,2}[t,\i)\Big\}\\
\ns\ds\qq\qq\subseteq\sU_{ad}[t,\i)\equiv\Big\{u(\cd)\in\sU^{\m,2}[t,\i)\bigm|J(t,x;u(\cd))
\hb{ is well-defined}\Big\}.\ea\ee
where $X(\cd)=X(\cd\,;t,x,(\Th,v(\cd)))$.
\ep

\it Proof. \rm Let $u(\cd)\in\sU^{\m,2}[t,\i)$, we have $X(\cd)\in\cX^{\m,2}[t,\i)$. Write
$$\ba{ll}
\ds dX(s)=[AX(s)+Bu(s)]ds+[CX(s)+Du(s)]dW(s)\\
\ns\ds\qq\q=[A^\Th X(s)+B(u(s)-\Th X(s))]ds+[C^\Th X(s)+D(u(s)-\Th X(s))]dW(s),\ea$$
with $\Th\in\hTh^\m[A,C;B,D]$. By letting
$$v(s)=u(s)-\Th X(s),\qq s\in[t,\i),$$
we see that
$(\Th,v(\cd))\in\hTh^\m[A,C;B,D]\times\sU^{\m,2}[t,\i)$. This gives the first inclusion. Next, let $(\Th,v(\cd))\in\hTh^\m[A,C;B,D]\times\sU^{\m,2}[t,\i)$. Then we write:
$$\ba{ll}
\ds dX(s)=[A^\Th X(s)+Bv(s)]ds+[C^\Th X(s)+Dv(s)]dW(s)\\
\ns\ds\qq\q=[AX(s)+B(\Th X(s)+v(s))]ds+[CX(s)+D(\Th X(s)+v(s))]dW(s).\ea$$
Then by setting
$$u(s)=\Th X(s)+v(s),\qq s\in[t,\i),$$
we see that the second inclusion hold. \endpf

\ms

We may further assume the following.

\ms

{\bf(H4)} Let $A,C\in\dbR^{n\times n}$, $B,D\in\dbR^{n\times m}$. Let system $[A,C;B,D]$ be weighted $L^2$-stabilizable with
$$0\in\hTh^\m[A,C;B,D].$$
Moreover, $Q\in\dbS_{++}^n$ and $R\in\dbS_{++}^m$ such that
\bel{Q>0*}Q-S^\top R^{-1}S\in\dbS_{++}^n,\ee
and $q(\cd)$ and $r(\cd)$ are deterministic with $|q(\cd)|+|r(\cd)|\in L^{\m,2}(0,\i;\dbR)$.

\ms

Under (H4), we see that $[A,C]$ is weighted $L^2$-stable.

\bp{} \sl Let {\rm(H4)} hold. Then
\bel{U=U}\sU_{ad}[t,\i)=\sU^{\m,2}[t,\i),\qq t\in[0,\i).\ee
\ep

\it Proof. \rm From the above proposition, we see that suffices to prove
$$\sU_{ad}[t,\i)\subseteq\sU^{\m,2}[t,\i).$$
For $(t,x)\in[0,\i)\times\dbR^n$ and $u(\cd)\in\sU_{a}[t,\i)$, we may let
$$X(s)\equiv X(s;t,x,u(\cd))=\G^0(s,t)x+\G(s,t)u(\cd),\qq s\in[t,\i),$$
Note that
$$\G^0(\cd\,,t)x,\q\G(\cd\,,t)u(\cd)\in L^{\m,2}_\dbF(t,\i;\dbR^n),\qq\forall (x,u(\cd))\in\dbR^n\times\sU^{\m,2}[t,\i).$$
Then one has
\bel{J=F}\ba{ll}
\ns\ds J(t,x;u(\cd))=\dbE\int_t^\i{\m(\t)\over\m(t)}\(\lan Q[\G^0(\t,t)x+\G(\t,t)u(\cd)],\G^0(\t,t)x+\G(\t,t)u(\cd)\ran\\
\ns\ds\qq\qq\qq\qq\qq+2\lan S[\G^0(\t,t)x+\G(\t,t)u(\cd)],u(\t)\ran+\lan R u(\t),u(\t)\ran\\
\ns\ds\qq\qq\qq\qq\qq+2\lan q(\t),\G^0(\t,t)x+\G(\t,t)u(\cd)\ran+2\lan r(\t), u(\t)\ran\)d\t\\
\ns\ds\qq\qq=\dbE\int_t^\i{\m(\t)\over\m(t)}\(\lan[R+\G(\t,t)^*Q\G(\t,t)+S\G(\t,t)
+\G(\t,t)^*S^\top]u(\cd),u(\cd)\ran\\
\ns\ds\qq\qq\qq\qq\qq+2\lan[\G(\t,t)^*Q\G^0(\t,t)+\G^0(\t,t)]x+\G(\t,t)^*
q(\t)+r(\t),u(\cd)\ran\\
\ns\ds\qq\qq\qq\qq\qq+\lan\G^0(\t,t)^*Q\G^0(\t,t)x,x\ran
+2\lan\G^0(\t,t)^*q(\t),x\ran\)d\t\\
\ns\ds\qq\qq\equiv\lan\F_2u(\cd),u(\cd)\ran+2\lan\F_1,u(\cd)\ran+\F_0,\ea\ee
where
$$\ba{ll}
\ns\ds\lan\F_2u(\cd),u(\cd)\ran=\dbE\int_t^\i{\m(\t)\over\m(t)}\lan[R+\G(\t,t)^* Q\G(\t,t)+S\G(\t,t)+\G(\t,t)^*S^\top]u(\cd),u(\cd)\ran d\t,\\
\ns\ds\lan\F_1,u(\cd)\ran=\dbE\int_t^\i{\m(\t)\over\m(t)}\lan[\G(\t,t)^*Q\G^0
(\t,t)+\G^0(\t,t)]x+\G(\t,t)^*q(\t)+r(\t),u(\t)\ran d\t,\\
\ns\ds\F_0=\dbE\int_t^\i{\m(\t)\over\m(t)}\(\lan\G^0(\t,t)^*Q\G^0(\t,t)x,x\ran
+2\lan\G^0(\t,t)^*q(\t),x\ran\)d\t.\ea$$
Thus, by (H4), we have some small $\e>0$,
\bel{F_2}\ba{ll}
\ns\ds\lan\F_2u(\cd),u(\cd)\ran=\dbE\[\int_t^\i{\m(\t)\over\m(t)}\(\e|u(\t)|^2+|(R-\e I)^{1\over2}[u(\t)
+(R-\e I)^{-1}S\G_\t u(\cd)|^2\\
\ns\ds\qq\qq\qq\qq\qq\qq+\lan(Q-S^\top(R-\e I)^{-1}S)\G_\t u(\cd),\G_\t u(\cd)\ran\)d\t\\
\ns\ds\qq\qq\qq\ges\e\dbE\int_t^\i{\m(\t)\over\m(t)}|u(\t)|^2d\t=\e\|u\|^2_{L^2(\O;
L^{\m,2}_\dbF(t,\i;\dbR^m))}.\ea\ee
Then we are done.\endpf

Now, we are ready to introduce the following:

\ms

\bf Problem (LQ). \rm Let (H4) hold. For any initial pair $(t,x)\in[0,\i)\times\dbR^n$, find $\bar u(\cd)\in\sU^{\m,2}[t,\i)$ such that
\bel{closed-sol}J(t,x;\bar u(\cd))=\inf_{u(\cd)\in\sU^{\m,2}[t,\i)}J(t,x;u(\cd))\equiv V(t,x).\ee
Any $\bar u(\cd)$ satisfying the above is called an {\it open-loop optimal control}. When this happens, Problem (LQ) is said to be {\it open-loop solvable}.

\ms

Next, for any $t\in[0,\i)$, and feedback strategy $(\Th,v(\cd))\in\hTh^\m[A,C;B,D]\times\sU^{\m,2}[t,\i)$, and any $x\in\dbR^n$, let $X(\cd)=X(\cd\,;t,x,(\Th,v(\cd)))$ be the solution to \rf{closed}. Then the above problem has another version.

\bde{} Problem (LQ) is said to be closed-loop solvable at $t\in[0,\i)$ if there exists a feedback strategy $(\Th,v(\cd))\in\hTh^\m[A,C;B,D]\times\sU^{\m,2}[t,\i)$ such that
$$J(t,x;(\bar\Th,\bar v(\cd)))\les J(t,x;(\Th,v(\cd))),\qq x\in\dbR^n.$$

\ede

We have the following simple result.

\bp{} \sl Let {\rm(H4)} hold and let $(\bar\Th,\bar v(\cd))\in\hTh^\m[A,C;B,D]\times\sU^{\m,2}[t,\i)$ be a closed-loop optimal strategy at $t\in[0,\i)$. Then
$$J(t,x;(\bar\Th,\bar v(\cd)))\les J(t,x;u(\cd)),\qq u(\cd)\in\sU^{\m,2}[t,\i),\q x\in\dbR^n.$$

\ep

\it Proof. \rm For given $(t,x)\in[0,\i)\times\dbR^n$, for any $u(\cd)\in\sU^{\m,2}[t,\i)$, let $X(\cd)=X(\cd\,;t,x,u(\cd))$, there exists a $(\Th,v(\cd))\in\hTh^\m[A,C;B,D]\times\sU^{\m,2}[t,\i)$ such that
$$u(s)=\Th X(s)+v(s),\qq s\in[t,\i).$$
Then
$$J(t,x;(\bar\Th,\bar v(\cd)))\les J(t,x;(\Th,v(\cd)))=J(t,x;u(\cd)).$$
This proves our proposition. \endpf

\ms

\section{Equivalence}

Directly further approach to Problem (LQ) seems to be difficult. Thus, we approach it indirectly instead. In this section, we will present an equivalence between the optimal control with recursive cost functional and the classical LQ problems. Recall the definition of $\m(\cd)$ from \rf{m}. Then
\bel{en}\sqrt{\m(s)\over\m(t
)}=e^{-{2E+F^2\over4}(s-t)-{F\over2}[W(s)-W(t)]}
\equiv e^{\n(s,t)},\ee
with
\bel{n}\n(s,t)=-{2E+F^2\over4}(s-t)-{F\over2}[W(s)-W(t)].\ee
Thus,
\bel{eqn}\left\{\2n\ba{ll}
\ds d\n(s,t)=-{2E+F^2\over4}ds-{F\over2}dW(s),\qq s\ges t,\\
\ns\ds\n(t,t)=0.\ea\right.\ee
By It\^o's formula, we have
\bel{d(en)}\ba{ll}
\ns\ds d\(e^{\n(s,t)}\)=e^{\n(s,t)}\(-{2E+F^2\over4}ds-{F\over2}dW(s)\)+{1\over2}
e^{\n(s,t)}{F^2\over4}ds\\
\ns\ds\qq\qq=e^{\n(s,t)}\(-{4E+F^2\over8}ds-{F\over2}dW(s)\).\ea\ee
Now, we let
\bel{wtXu}\ba{ll}
\ns\ds\wt X(s)=\sqrt{\m(s)\over\m(t)}X(s)\equiv e^{\n(s,t)}X(s),\\
\ns\ds\wt u(s)=\sqrt{\m(s)\over\m(t)}u(s)\equiv e^{\n(s,t)}u(s),\ea\qq s\in[t,\i).\ee
Thus,
\bel{dwtX}\ba{ll}
\ns\ds d\wt X(s)=d\(e^{\n(s,t)}X(s)\)\\
\ns\ds=e^{\n(s,t)}X(s)\(-{4E\1n+\1n F^2\over8}ds-{F\over2}dW(s)\)
\1n+\1n e^{\n(s,t)}\([AX(s)\1n+\1n Bu(s)]ds\1n+\1n[CX((s)\1n+\1n Du(s)]dW(s)\)\\
\ns\ds\qq-{F\over2}e^{\n(s,t)}[CX(s)+Du(s)]ds\\
\ns\ds=\[\(A-{F\over2}C-{4E+F^2\over8}I\)\wt X(s)+\(B-{F\over2}D\)\wt u(s)\]ds+\[\(C-{F\over2}I\)\wt X(s)+D\wt u(s)\]dW(s)
\ea\ee
Hence, we end up with the following control system:
\bel{wtX}\left\{\2n\ba{ll}
\ns\ds d\wt X(s)=\[\wt A\wt X(s)+\wt B\wt u(s)\]ds+\[\wt C\wt X(s)+\wt D\wt u(s)\]dW(s),\qq s\in[t,\i),\\
\ns\ds\wt X(t)=x,\ea\right.\ee
where
\bel{wtA}\left\{\2n\ba{ll}
\ns\ds\wt A=A-{F\over2}C-{4E+F^2\over8}I,\qq\wt B=B-{F\over2}D,\\
\ns\ds\wt C=C-{F\over2}I,\qq\wt D=D,\ea\right.\ee
Thus, system \rf{wtX} should be denoted by $[\wt A,\wt C;\wt B,\wt D]$. Also, we see that the cost is
\bel{J=J}\ba{ll}
\ns\ds J(t,x;u(\cd))=\dbE\int_t^\i {\m(\t)\over\m(t)}\(\lan QX(\t),X(\t)\ran+2\lan SX(\t),u(\t)\ran+\lan Ru(\t),u(\t)\ran\\
\ns\ds\qq\qq\qq\qq\qq+2\lan q(\t),X(\t)\ran+2\lan r(\t),u(\t)\ran\)d\t\\
\ns\ds\qq=\dbE\int_t^\i\3n\(\lan Q\wt X(\t),\wt X(\t)\ran+2\lan S\wt X(\t),\wt u(\t)\ran+\lan R\wt u(\t),\wt u(\t)\ran+2\lan\wt q(\t),\wt X(\t)\ran+2\lan\wt r(\t),\wt u(\t)\ran\)d\t\\
\ns\ds\qq\equiv\wt J(t,x;\wt u(\cd)),\ea\ee
with
\bel{wtqr}\wt q(s)=\sqrt{\m(s)\over\m(t)}q(s)\equiv e^{\n(s,t)}q(s),\q\wt r(s)=\sqrt{\m(s)\over\m(t)}r(s)\equiv e^{\n(s,t)}r(s),\q s\ges t,\ee
which are random (although $q(\cd)$ and $r(\cd)$ are deterministic, which will be important below).
Consequently, system \rf{wtX} with cost $\wt J(t,x;\wt u(\cd))$ leads to a classical LQ stochastic optimal control in an infinite time horizon. Then we can pose the following:

\ms

\bf Problem ($\wt{\bf LQ}$). \rm Let (H4) hold. For any initial pair $(t,x)\in[0,\i)\times\dbR^n$, find $\bar{\wt u}(\cd)\in\sU^2[t,\i)$ such that
\bel{closed-sol}\wt J(t,x;\bar{\wt u}(\cd))=\inf_{\wt u(\cd)\in\sU^2[t,\i)}\wt J(t,x;\wt u(\cd))\equiv\wt V(t,x).\ee
Any $\bar{\wt u}(\cd)$ satisfying the above is an open-loop optimal control. When this happens, Problem ($\wt{\rm LQ}$) is said to be open-loop solvable.

\ms

Clearly, we see that Problem (LQ) is open-loop solvable if and only if so is Problem ($\wt{\rm LQ}$). Next, for any $t\in[0,\i)$, we can define the $L^2$-stability of $[\wt A,\wt C]$ and the $L^2$-stabilzabilty of $[\wt A,\wt C;\wt B,\wt D]$ as in \cite{Sun-Yong 2020}. We let $\hTh[\wt A,\wt C;\wt B,\wt D]$ be the set of all $L^2$-stabilizers of the controlled system $[\wt A,\wt C;\wt B,\wt D]$. In what follows, for any $(t,x)\in[0,\i)\times\dbR^n$ and state feedback strategy $(\wt\Th,\wt v(\cd))\in\hTh[\wt A,\wt C;\wt B,
\wt D]\times\sU^2[t,\i)$, we let $\wt X(\cd)=\wt X(\cd\,;t,x,(\wt\Th,\wt v(\cd)))$ be the solution to the following:
\bel{wtX*}\left\{\2n\ba{ll}
\ns\ds d\wt X(s)=\[(\wt A+\wt B\wt\Th)\wt X(s)+\wt B\wt v(s)\]ds+\[(\wt C+\wt D\wt\Th)\wt X(s)+\wt D\wt v(s)\]dW(s),\qq s\in[t,\i),\\
\ns\ds\wt X(t)=x.\ea\right.\ee
This system is referred to as the closed-loop system under strategy $(\wt\Th,\wt v(\cd))$ and the initial state $x\in\dbR^n$ at $t\in[0,\i)$.
\bde{} Problem ($\wt{\rm LQ}$) is said closed-loop solvable at $t\in[0,\i)$ if there exists a feedback strategy $(\bar{\wt\Th},\bar{\wt v}(\cd))\in\hTh[\wt A,\wt C;\wt B,\wt D]\times\sU^2[t,\i)$ such that
$$\wt J(t,x;(\bar{\wt\Th},\bar{\wt v}(\cd)))\les\wt J(t,x;(\wt\Th,\wt v(\cd))),\qq x\in\dbR^n.$$

\ede

Clearly, Problem (LQ) is closed-loop solvable if and only if so is Problem ($\wt {\rm LQ}$). We now present the following result.

\bp{Prop.4.2.} \sl Let {\rm(H4)} hold.

\ms

{\rm(i)} System $[A,C]$ is weighted $L^2$-stable if and only if $[\wt A,\wt C]$ is $L^2$-stable.

\ms

{\rm(ii)} System $[A,C;B,D]$ is weighted $L^2$-stabilizable if and only if $[\wt A,\wt C;\wt B,\wt D]$ is $L^2$-stabilizable.

\ep

\it Proof. \rm (i) system $[\wt A,\wt C]$ is $L^2$-stable if and only if for any $\L\in\dbS_{++}^n$, the following Lyapunov equation admits a solution $P\in\dbS^n_{++}$:
$$\ba{ll}
\ns\ds0=P\wt A+\wt A^\top P+\wt C^\top P\wt C+\L\\
\ns\ds\q=P\(A-{F\over2}C-{4E+F^2\over8}I\)+\(A-{F\over2}C-{4E+F^2\over8}I\)^\top P+\(C-{F\over2}I\)^\top P\(C-{F\over2}I\)+\L\\
\ns\ds\q=PA+A^\top P+C^\top PC-{F\over2}(PC+C^\top P)-{4E+F^2\over4}P-{F\over2}(C^\top P+PC)+{F^2\over4}P+\L\\
\ns\ds\q=-EP-F(PC+C^\top P)+PA+A^\top P+C^\top PC+\L.\ea$$
This is true if and only if $[A,C]$ is weighted $L^2$-stable.

\ms

(ii) Controlled system $[\wt A,\wt C;\wt B,\wt D]$ is $L^2$-stabilizable if and only for some $\Th\in\dbR^{n\times m}$, and any $\L\in\dbS^n_{++}$, the following admits a solution $P\in\dbS^n_{++}$:
$$\ba{ll}
\ns\ds0=P(\wt A+\wt B\Th)+(\wt A+\wt B\Th)^\top P+(\wt C+\wt D\Th)^\top P(\wt C+\wt D\Th)+\L\\
\ns\ds\q=P\[\(A-{F\over2}C-{4E+F^2\over8}I\)+\(B-{F\over2}D\)\Th\]
+\[\(A-{F\over2}C-{4E+F^2\over8}I\)+\(B-{F\over2}D\)\Th\]^\top P\\
\ns\ds\qq+\[\(C-{F\over2}I\)+D\Th\]^\top P\[\(C-{F\over2}I\)+D\Th\]+\L\\
\ns\ds\q=P\[A+B\Th-{F\over2}(C-D\Th)-{4E+F^2\over8}I\]
+\[A+B\Th-{F\over2}(C-D\Th)-{4E+F^2\over8}I\]^\top P\\
\ns\ds\qq+\[(C-D\Th)-{F\over2}I\]^\top P\[(C-D\Th)-{F\over2}I\]+\L\\
\ns\ds\q=-EP-F[P(C+D\Th)+(C+D\Th)^\top P]+P(A+B\Th)+(A+B\Th)^\top P+(C+D\Th)^\top P(C+D\Th)+\L.\ea$$
This is true if and only if $[A,C;B,D]$ is weighted $L^2$-stabilizable. \endpf

\section{A Translation}

From the results from the previous section, to get the results for Problem (LQ), it suffices to translate those for Problem ($\wt{\rm LQ}$). In this section, we are going to carry out this translation. For Problem ($\wt{\rm LQ}$), it is a classical one and we have the following result (see \cite{Sun-Yong 2020} p.90, for details).

\bt{wt} \sl Let {\rm(H4)} hold. $\wt Q=Q, \wt S=S, \wt R=R.$

\ms

{\rm(i)} Let $(t,x)\in[0,\i)\times\dbR^n$ be given. Let $(\wt X(\cd),\wt u(\cd))$ be a state-control pair of \rf{wtX} corresponding the initial pair $(t,x)$. Then this pair is an open-loop optimal pair of Problem {\rm($\wt{\rm LQ}$)} if and only if the following FBSDE
\bel{FBSDE}\left\{\2n\ba{ll}
\ns\ds d\wt X(s)=[\wt A\wt X(s)+\wt B\wt u(s)]ds+[\wt C\wt X(s)+\wt D\wt u(s)]dW(s),\\
\ns\ds d\wt Y(s)=-[\wt A^\top\wt Y(s)+\wt C^\top\wt Z(s)+\wt Q\wt X(s)+\wt S^\top\wt u(s)+\wt q(s)]ds+\wt Z(s)dW(s),\\
\ns\ds\wt X(t)=x,\qq\wt Y(\cd)\in L^2_{\dbF^t}(t,\i;\dbR^n).\ea\right.\ee
admits an adapted solution $(\wt X(\cd),\wt u(\cd),\wt Y(\cd),\wt Z(\cd))$, with the stationary condition:
\bel{stationary}\wt B^\top\wt Y(s)+\wt D^\top\wt Z(s)+\wt S\wt X(s)+\wt R\wt u(s)+\wt r(s)=0,\qq s\in[t,\i),\as\ee

\ms

{\rm(ii)} Problem {\rm($\wt{\rm LQ}$)} is closed-loop solvable with the closed-loop optimal strategy $(\bar{\wt\Th},\bar{\wt v}(\cd))$ if and only if the following ARE:
\bel{ARE}\ba{ll|}
\ns\ds\wt P\wt A+\wt A^\top\wt P+\wt C^\top\wt P\wt C+\wt Q\\
\ns\ds-(\wt B^\top\wt P+\wt D^\top\wt P\wt C+\wt S)^\top (\wt R+\wt D^\top\wt P\wt D)^{-1}(\wt B^\top\wt P
+\wt D^\top\wt P\wt C+\wt S)=0,\ea\ee
admits a unique $L^2$-stabilizing solution $\wt P\in\dbS_{++}^n$, with respect to the system $[\wt A,\wt C;\wt B,\wt D]$, i.e.,
\bel{Th}\bar{\wt\Th}=-(\wt R+\wt D^\top\wt P\wt D)^{-1}(\wt B^\top\wt P+\wt D^\top\wt P\wt C+\wt S)\in\hTh[\wt A,\wt C;\wt B,\wt D],\ee
and
\bel{wteta}\bar{\wt v}(s)=-(\wt R+\wt D^\top\wt P\wt D)^{-1}[\wt B^\top\wt\eta(s)+D^\top\wt\z(s)+\wt r(s)],\ee
with $(\wt\eta(\cd),\wt\z(\cd))\in L^2_{\dbF^t}(\O;L^2(t,\i;\dbR))\times L^2_{\dbF^t}(\O;L^2(t,\i;\dbR))$ solving
\bel{eta}d\wt\eta(s)=-[(\wt A+\wt B\bar{\wt\Th})^\top\wt\eta(s)+(\wt C+\wt D\bar{\wt\Th})^\top\wt\z(s)+\bar{\wt\Th}^\top\wt r(s)+\wt q(s)]ds+\wt\z(s)dW(s), \qq s\in[t,\i),\ee
with
$$\dbE\[\sup_{s\in[t,\i)}|\wt\eta(s)|^2+\(\int_t^\i|\wt\z(s)|^2ds\)\]\les K\dbE\int_t^\i\(|\wt q(s)|^2+|\wt r(s)|^2\)ds.$$

\et

Note in the above that although $q(\cd)$ and $r(\cd)$ are assumed to be deterministic, $\wt q(\cd)$ and $\wt r(\cd)$ are random. Therefore, in the above, $\wt\z(\cd)$ might exist and could be nonzero. Note the system $[\wt A+\wt B\bar{\wt\Th},\wt C+\wt D\bar{\wt\Th}]$ is $L^2$-stable. Thus, $[A+B\bar{\wt\Th}, C+D\bar{\wt\Th}]$ is weighted $L^2$-stable.

\ms

We now translate the above in terms of control system $[A,C;B,D]$  and the recursive functional.

To this end, we first write FBSDE \rf{FBSDE} and stationary condition \rf{stationary} as follows. Let
\bel{wtn}\wt\n(s,t)=-\n(s,t)={2E+F^2\over4}(s-t)+{F\over2}[W(s)-W(t)],\qq s\in[t,\i).\ee
Then,
\bel{eq-n}\left\{\2n\ba{ll}
\ds d\wt\n(s,t)={2E+F^2\over4}ds+{F\over2}dW(s),\qq s\ges t,\\
\ns\ds\wt\n(t,t)=0.\ea\right.\ee
Consequently,
\bel{d(ewtn)}\ba{ll}
\ns\ds d\(e^{\wt\n(s,t)}\)=e^{\wt\n(s,t)}\({2E+F^2\over4}ds+{F\over2}dW(s)\)+{1\over2}
e^{\wt\n(s,t)}{F^2\over4}ds\\
\ns\ds\qq\qq=e^{\wt\n(s,t)}\({4E+3F^2\over8}ds+{F\over2}dW(s)\).\ea\ee
From \rf{wtX}, we have
\bel{Xu*}\ba{ll}
\ns\ds X(s)=e^{\wt\n(s,t)}\wt X(s)\equiv\sqrt{\m(t)\over\m(s)}\wt X(s),\\
\ns\ds u(s)=e^{\wt\n(s,t)}\wt u(s)\equiv\sqrt{\m(t)\over\m(s)}\wt u(s),\ea\qq s\ges t.\ee
Also, we define
\bel{YZ}\ba{ll}
\ns\ds Y(s)=e^{\wt\n(s,t)}\wt Y(s)\equiv\sqrt{\m(t)\over\m(s)}\wt Y(s),\\
\ns\ds Z(s)=e^{\wt\n(s,t)}\(\wt Z(s)+{F\over2}\wt Y(s)\)\equiv\sqrt{\m(t)\over\m(s)}\(\wt Z(s)+{F\over2}\wt Y(s)\).\ea\qq s\ges t.\ee
Hence,
$$\ba{ll}
\ns\ds dX(s)=d\(e^{\wt\n(s,t)}\wt X(s)\)=\({4E+3F^2\over8}ds+{F\over2}dW(s)\)X(s)+\(\wt AX(s)+\wt Bu(s)\)ds\\
\ns\ds\qq\qq+\(\wt CX(s)+\wt Du(s)\)dW(s)+{F\over2}\(\wt CX(s)+\wt Du(s)\)ds\\
\ns\ds\qq=\[\(\wt A+{4E+3F^2\over8}+{F\over2}\wt C\)X(s)+\(\wt B+{F\over2}\wt D\)u(s)\]ds+\[\(\wt C+{F\over2}I\)X(s)+Du(s)\]dW(s)\\
\ns\ds\qq=[AX(s)+Bu(s)]ds+[CX(s)+Du(s)]dW(s),\ea$$
which is expected. At the same time, we have
$$\ba{ll}
\ns\ds dY(s)=d\(e^{\wt\n(s,t)}\wt Y(s)\)=\({4E+3F^2\over8}ds+{F\over2}dW(s)\)Y(s)\\
\ns\ds\qq\qq-\(\wt A^\top Y(s)+\wt C^\top[Z(s)-{F\over2}Y(s)]+QX(s)+\wt S^\top u(s)\)ds\\
\ns\ds\qq\qq+\[Z(s)-{F\over2}Y(s)\]dW(s)+{F\over2}\[Z(s)-{F\over2}Y(s)+(C-FI)^\top Z(s)+QX(s)+S^\top u(s)\]ds\\
\ns\ds\qq=-\[\(\wt A-{4E+3F^2\over8}I-{F\over2}\wt C+{F^2\over4}I\)^\top Y(s)+\(\wt C-{F\over2}\)^\top Z(s)+QX(s)+\wt S^\top u(s)\)\]ds\\
\ns\ds\qq\qq+Z(s)dW(s)\\
%
%
\ns\ds\qq=-\[\(A-EI-FC\)^\top Y(s)+(C-FI)^\top
Z(s)+QX(s)+S^\top u(s)\]ds\\
\ns\ds\qq\qq+Z(s)dW(s)\ea$$
On the other hand, we see that, by (5.11)
$$\dbE \int_t^\i {\m(s)\over \m(t)}|Y(s)|^2 ds = \dbE \int_t^\i |\wt Y(s)|^2 ds<\i. $$
Hence, we must have
$$|Y(\cd)|\in L^2_{\dbF^t}(\O;L^{\m,2}(t,\i;\dbR)), \quad \lim_{t\to \i}Y(t)=\lim_{t\to \i}\wt Y(t)=0, \quad a.s.$$
The stationary condition reads
\bel{stationary*}\ba{ll}
\ns\ds0=\wt B^\top Y(s)+D^\top\(Z(s)-{F\over2}Y(s)\)+SX(s)+Ru(s)\\
\ns\ds\q=(B-FD)^\top Y(s)+D^\top Z(s)+SX(s)+Ru(s)\ea\ee
Hence, we have

\bt{open} \sl Let {\rm(H4)} hold and $(X(\cd),u(\cd))$ be state-control pair of \rf{state}. Then this pair is an open-loop optimal pair of Problem {\rm(LQ)} if and only if the following FBSDE
\bel{FBSDE*}\left\{\2n\ba{ll}
\ds dX(s)=[AX(s)+Bu(s)]ds+[CX(s)+Du(s)]dW(s),\\
\ns\ds dY(s)=-[(A-EI-FC)^\top Y(s)+(C-FI)^\top Z(s)\\
\ns\ds\qq\qq\qq\qq+QX(s)+S^\top u(s)+q(s)]ds+Z(s)dW(s),\qq s\in[t,\i),\\
\ns\ds X(t)=x,\q|Y(\cd)|\in L^2_\dbF(\O;L^{\m,2}(t,\i;\dbR)),\ea\right.\ee
admits an adapted solution $(X(\cd),u(\cd),Y(\cd),Z(\cd))$ with the stationary condition:
\bel{stationary*}(B-FD)^\top Y(s)+D^\top Z(s)+SX(s)+Ru(s)+r(s)=0,\qq s\in[t,T].\ee

\et

Next, for closed-loop solvability of Problem (LQ), we first note that the ARE \rf{ARE} can be written as
\bel{ARE*}\ba{ll|}
\ns\ds0=\wt P\(A-{F\over2}C-{4E+F^2\over8}I\)+\(A-{F\over2}C-{4E+F^2\over8}I\)^\top\wt P+\(C-{F\over2}I\)^\top\wt P\(C-{F\over2}I\)+Q\\
\ns\ds\qq\qq-\[\(B-{F\over2}D\)^\top\wt P+D^\top\wt P\(C-{F\over2}I\)+S\]^\top(R+D^\top\wt PD)^{-1}\\
\ns\ds\qq\qq\q\cd\[\(B-{F\over2}D\)^\top\wt P+D^\top\wt P\(C-{F\over2}I\)+S\]\\
\ns\ds\q=\wt PA+A^\top\wt P+C^\top\wt PC+Q-E\wt P-F(\wt PC+C^\top\wt P)\\
\ns\ds\qq-(B^\top\wt P+D^\top\wt PC+S-FD^\top\wt P)^\top(R+D^\top\wt PD)^{-1}(B^\top\wt P+D^\top\wt PC+S-FD^\top\wt P).\ea\ee
Thus, this ARE has an $L^2$-stabilizing solution $\wt P$ (with respect to the system $[\wt A,\wt C;\wt B,\wt D]$), meaning that
$$\bar{\wt\Th}=-(\wt R+\wt D^\top\wt P\wt D)^{-1}\(\wt B^\top\wt P+\wt D^\top\wt P\wt C+\wt S\)\in
\hTh[\wt A,\wt C;\wt B,\wt D].$$
This is true if and only if for any $\L\in\dbS_{++}^n$, the Lyapunov equation %
$$\ba{ll}
\ns\ds P(\wt A+\wt B\bar{\wt\Th})+(\wt A+\wt B\bar{\wt\Th})^\top P+(\wt C+D\bar{\wt\Th})^\top P(\wt C+D\bar{\wt\Th})+\L=0,\ea$$
admits a solution $P\in\dbS_{++}^n$. Then
$$\ba{ll}
\ns\ds P\[A-{F\over2}C-{4E+F^2\over8}I+\(B-{F\over2}D\)\bar{\wt\Th}\]
+\[A-{F\over2}C-{4E+F^2\over8}I+\(B-{F\over2}D\)\bar{\wt\Th}\]^\top P\\
\ns\ds\qq+\[\(C-{F\over2}I\)+D\bar{\wt\Th}\]^\top P\[\(C-{F\over2}I\)+D\bar{\wt\Th}\]+\L=0,\ea$$
This gives
$$\ba{ll}
\ns\ds0=P(A+B\bar{\wt\Th})+(A+B\bar{\wt\Th})^\top P+(C+D\bar{\wt\Th})^\top P(C+D\bar{\wt\Th})+\L\\
\ns\ds\qq-{F\over2}PC-{4E+F^2\over8}P-{F\over2}PD\bar{\wt\Th}-{F\over2}C^\top P
-{4E+F^2\over8}P-{F\over2}\bar{\wt\Th}^\top D^\top P\\
\ns\ds\qq-{F\over2}P(C+D\bar{\wt\Th})+{F\over2}(C+D\bar{\wt\Th})^\top P+{F^2\over 4}P\\
\ns\ds\q=-EP+P(A+B\wt\Th)+(A+B\bar{\wt\Th})^\top P+(C+D\bar{\wt\Th})^\top P(C+D\bar{\wt\Th})+\L\\
\ns\ds\qq-F[P(C+D\bar{\wt\Th})+(C+D\bar{\wt\Th})^\top P],\ea$$
which shows that $\bar{\wt\Th}\in\hTh^\m[A,C;B,D]$ (This also follows from Proposition \ref{Prop.4.2.}). We also have
\bel{wtTh}\ba{ll}
\ns\ds\bar{\wt\Th}=-(\wt R+\wt D^\top\wt P\wt D)^{-1}\(\wt B^\top\wt P+\wt D^\top\wt P\wt C+\wt S\)\\
\ns\ds\q=-(R+D^\top\wt PD)^{-1}\((B^\top-{F\over2}D^\top)\wt P+D^\top\wt P(C-{F\over2}I)+S\)\\
\ns\ds\q=-(R+D^\top\wt PD)^{-1}(B^\top\wt P+D^\top\wt PC+S-FD^\top\wt P).\ea\ee
Now, we define
$$\eta(s)=e^{\wt\n(s,t)}\wt\eta(s),\q\z(s)=e^{\wt\n(s,t)}\wt\z(s),\qq s\in[t,\i).$$
Then, $\eta(\cd)\in\cX^{\m,2}[t,\i)$ and $\z(\cd)\in L_{\dbF^t}^2(\O;L^{\m,2}(t,\i:\dbR^n)$. By It\^o's formula,
$$\ba{ll}
\ns\ds d\eta(s)=d\(e^{\wt\n(s,t)}\wt\eta(s)\)\\
\ns\ds\qq=\({4E+3F^2\over8}ds+{F\over2}dW(s)\)\eta(s)-\[\(\wt A+\wt B\bar{\wt\Th}\)^\top\eta(s)+\(\wt C+\wt D\bar{\wt\Th}\)^\top\z(s)+\bar{\wt\Th}^\top r(s)+q(s)\]ds\\
\ns\ds\qq\qq+\z(s)dW(s)+{F\over2}\z(s)ds\\
\ns\ds\qq=-\[\(\wt A-{4E+3F^2\over8}I+\wt B\bar{\wt\Th}\)^\top\eta(s)
+\(\wt C+D\bar{\wt\Th}-{F\over2}I\)\z(s)+\bar{\wt\Th}^\top r(s)+q(s)\]ds\\
\ns\ds\qq\qq+\({F\over2}\eta(s)+\z(s)\)dW(s).\ea$$
Since all the coefficients are deterministic, we must have
$${F\over2}\eta(s)+\z(s)=0.$$
Hence,
\bel{eta*}\ba{ll}
\ns\ds d\eta(s)=-\[\(\wt A-{4E+3F^2\over8}I+\wt B\bar{\wt\Th}\)^\top\eta(s)
+\(\wt C+D\bar{\wt\Th}-{F\over2}I\)^\top\(-{F\over2}\eta(s)\)+\bar{\wt\Th}^\top r(s)+q(s)\]ds\\
\ns\ds\qq\q=-\Big\{\[\(A-{F\over2}C-{4E+F^2\over8}I-{4E+3F^2\over8}I
+(B-{F\over2}D)\bar{\wt\Th}\)^\top\\
\ns\ds\qq\qq\q-\((C-{F\over2}I)+D\bar{\wt\Th}-{F\over2}I\)^\top{F\over2}\]\eta(s)
+\bar{\wt\Th}^\top r(s)+q(s)\Big\}ds\\
\ns\ds\qq\q=-\[\(A-EI-FC+(B-FD)\bar{\wt\Th}\)^\top\eta(s)
+\bar{\wt\Th}^\top r(s)+q(s)\]ds\ea\ee
Then by \rf{wteta},
$$\bar{\wt v}(s)=-(R+D^\top\wt PD)^{-1}\[\(B-{F\over2}D\)^\top\wt\eta(s)+D^\top\wt\z(s)+\wt r(s)\].$$
Thus,
\bel{barv}\ba{ll}
\ns\ds\bar v(s)=\sqrt{\m(t)\over\m(s)}\bar{\wt v}(s)=-(R+D^\top\wt PD)^{-1}\[\(B-{F\over2}D\)^\top\eta(s)+D^\top\z(s)+r(s)\]\\
\ns\ds\qq=-(R+D^\top\wt PD)^{-1}\[(B-FD)^\top\eta(s)+r(s)\],\q s\in[t,\i).\ea\ee
From the above, we have the following result.

\bt{closed} \sl Let {\rm(H4)} hold. Then Problem {\rm(LQ)} is closed-loop solvable with $(\bar{\wt\Th},\bar v(\cd))$ being a closed-loop optimal pair if and only if the following ARE
\bel{ARE0}\ba{ll}
\ns\ds\wt PA+A^\top\wt P+C^\top\wt PC+Q-E\wt P-F(\wt PC+C^\top\wt P)\\
\ns\ds\qq-(B^\top\wt P+D^\top\wt PC+S-FD^\top\wt P)^\top(R+D^\top\wt PD)^{-1}(B^\top\wt P+D^\top\wt PC+S-FD^\top\wt P)=0,\ea\ee
admits a weighted $L^2$-stabilizing solution $\wt P\in\dbS_{++}^n$, i.e.,
$\bar{\wt\Th}$ defined by \rf{wtTh} is in $\hTh^\m[A,C;B,D]$, $\bar v(\cd)$ defined by \rf{barv} is in $L^2_{\dbF^t}(\O;L^{\m,2}(t,\i;\dbR^m))$ and with $\eta(\cd)\in L^2_{\dbF^t}(\O;L^{\m,2}(t,\i;\dbR^n))$ being the solution of \rf{eta*}.

\et

\br{} It is pity for us that we are not able to directly prove that the equation \rf{eta*} is solvable by $\eta(\cd)$ (not via $\wt\eta(\cd)$).
\er

\section{An Extension}

In this section, we consider the following nonhomogeneous controlled system:
\bel{state2}\left\{\2n\ba{ll}
\ds dX(s)=[AX(s)+Bu(s)+b(s)]ds+[CX(s)+Du(s)+\si(s)]dW(s),\q s\in[t,T],\\
\ns\ds X(t)=x.\ea\right.\ee
The difference of the above with \rf{state} is $b,\si:[0,\i)\to\dbR^n$ appear. Note that they are assumed to be deterministic. The reason is the same as that for $q(\cd)$ and $r(\cd)$. The state process $X(\cd\,;t,x,u(\cd))$ can be decomposed
$$X(\cd\,;t,x,u(\cd))=X_0(\cd\,;t,x,u(\cd))+X(\cd\,;t,0,0)\equiv X_0(\cd)+\h X(\cd).$$
In the above, the two processes are those of state \rf{state2} with $X_0(\cd)$ corresponding to $b(\cd),\si(\cd)=0$ and $\h X(\cd)$ corresponding to $x,u(\cd)=0$.
Then
$$\ba{ll}
\ns\ds J(t,x;u(\cd))=\dbE\int_t^\i{\m(\t)\over\m(t)}f(\t)d\t\\
\ns\ds\qq\qq=\dbE\int_t^\i{\m(\t)\over\m(t)}\(\lan Q[X_0(\t)+\h X(\t)],X_0(\t)+\h X(\t)\ran+2\lan S[X_0(\t)+\h X(\t)],u(\t)\ran\\
\ns\ds\qq\qq\qq\qq\qq\qq\qq+\lan Ru(\t),u(\t)\ran+2\lan q(\t),X_0(\t)+\h X(\t)\ran+2\lan r(\t),u(\t)\ran\)d\t\\
\ns\ds\qq\qq=\dbE\int_t^\i{\m(\t)\over\m(t)}\(\lan QX_0(\t),X_0(\t)\ran+2\lan S X_0(\t),u(\t)\ran+\lan Ru(\t),u(\t)\ran\\
\ns\ds\qq\qq\qq\qq\qq\qq\qq+2\lan q(\t)+Q\h X(\t),X_0(\t)\ran+2\lan r(\t)+S\h X(\t),u(\t)\ran\\
\ns\ds\qq\qq\qq\qq\qq\qq\qq+\lan Q\h X(\t),\h X(\t)\ran+2\lan q(\t),\h X(\t)\ran\)d\t\\
\ns\ds\equiv J^0(t,x;u(\cd))+\f(t),\ea$$
where
\bel{J_0}\ba{ll}
\ns\ds J^0(t,x;u(\cd))=\dbE\int_t^\i{\m(\t)\over\m(t)}\(\lan QX_0(\t),X_0(\t)\ran+2\lan SX_0(\t),u(\t)\ran+\lan Ru(\t),u(\t)\ran\\
\ns\ds\qq\qq\qq\qq\qq\qq\qq+2\lan\h q(\t),X_0(\t)\ran+2\lan\h r(\t),u(\t)\ran\)d\t,\\
\ns\ds\f(t)=\dbE\int_t^\i{\m(\t)\over\m(t)}\(\lan Q\h X(\t),\h X(\t)\ran+2\lan q(\t),\h X(\t)\ran\)d\t,\ea\ee
with
$$\h q(s)=q(s)+Q\h X(s),\qq\h r(s)=r(s)+S\h X(s),\qq s\in[t,\i).$$
Since $\f(t)$ is independent of $(x,u(\cd))\in\dbR^n\times\sU^2[t,\i)$, minimizing the map $u(\cd)\mapsto J(t,x;u(\cd))$ is the same as minimizing the
map $u(\cd)\mapsto J^0(t,x;u(\cd))$. Hence, we may assume that
\bel{bsi}b(\cd)=\si(\cd)=0.\ee
The problem is reduced the one we discussed earlier.

\section{Conclusions}

We formulated a stochastic LQ problem with recurve cost functional in an infinite horizon. It is naturally involved BSDE in $L^1$, which leads to the notion of weighted $L^2$-stabilizabiity of the controlled system. This further leads to an equivalence of the considered problem with a classical one. Consequently, all the results for the classical LQ problem can be translated.

\ms


\begin{thebibliography}{99}

\bibitem{Briand-Hu 2005} P.~Briand and Y.~Hu, \it BSDE with quadratic growth and unbounded terminal value, \rm arXiv:math.PR/0504002 v1 1 Apr 2005.

\bibitem{Duffie-Epstein 1992a} D.~Duffie and L.~G.~Epstein, \it Stochastic differential utility, \sl Econometrica,  \rm 60 (1992), no. 2, 353-394.

\bibitem{Duffie-Epstein 1992b} D.~Duffie and L.~G.~Epstein, \it Asset pricing with stochastic differential utility, \sl Review Financial Studies, \rm 5 (1992), 411-436.

\bibitem{El Karoui-Peng-Quenez 1997} N.~El Karoui, S.~Peng, and M.~C.~Quenez, \it Backward stochastic differential equations in finance, \sl Math. Finance, \rm 7 (1997), 1--71.

\bibitem{Fan 2016} S.~Fan, \it Bounded solutions, $L^p$ ($p>1$) solutions and $L^1$ solutions for one dimensional BSDEs under general assumptions, \sl Stoch. Proc. Appl., \rm 126 (2016), 1511--1552.

\bibitem{Fan-Liu 2010} S.~Fan and D.~Liu, \it A class of BSDE with integrable parameters, \sl Stat. Probab. Letters, \rm 80 (2010), 2024--2031.

\bibitem{Fuhrman-Tessitore 2004} M.~Fuhrman and G.~Tessitore, \it Infinite horizon backward stochastic differential equations and elliptic equations in Hilbert spaces, \sl Ann. Probab.,
    \rm 32 (2004), 607--660.

\bibitem{Klimsaik-Rzymowski 2024} T.~Klimsiak and M.~Rzymowski, \it A priori estimates for multidimensional BSDEs with integrable data, \sl Electron. Commun. Probab., \rm 29 (2024), no.34.


\bibitem{Lazrak 2004} A.~Lazrak, \it Generalized stochastic differential utility and preference for information, \sl Ann. Appl. Probab., \rm 14 (2004), 2149--2175.

\bibitem{Lazrak-Quenez 2003} A.~Lazrak and M.~C.~Quenez, \it A generalized stochastic differential utility, \sl Math. Oper. Res., \rm 28 (2003), 154--180.

\bibitem{Li-Yong 2026} L.~Li and J.~Yong, \it Stochastic optimal linear quadratic controls with a recursive dost functional, \sl Evol. Equa Control Th., \rm to appear. arXiv preprint arXiv:2601.21748, 2026.

\bibitem{Luo-Li-Wei 2025} S.~Luo, X.~Li, and Q.~Wei, \it Infinite time horizon stochastic recursive control problems with jumps: dynamic programming and stochastic verification theorems, \sl. SIAM J. Control Optim.,
    \rm 63 (2025) 796--821.

\bibitem{Ma-Protter-Yong 1994} J.~Ma, P.~Protter, and J.~Yong, \it Solving forward-backward stochastic differential equations explicitly --- a four step scheme., \sl Probab. Theory Related Fields, \rm  98 (1994), 339--359.

\bibitem{Ma-Yong 1999} J.~Ma and J.~Yong, \sl Forward-Backward Stochastic Differential Equations and Their Applications, \rm Lecture Notes in Math., Vol. 1702, Springer-Verlag, 1999.

\bibitem{Mou-Yong 2006} L.~Mou and J.~Yong, \it Two-person zero-sum linear quadratic stochasrtic differential games by a Hilbert space method, \sl J. Ind. Manag. Optim., \rm 2 (2006), 95--117.

\bibitem{Peng 1992} S.~Peng, \it A generalized dynamic programming principle and Hamilton-Jacobi-Bellman equation, \sl Stoch. Stoch. Rep., \rm 38 (1992), 119--134.

\bibitem{Peng 1997} S.~Peng, \it Backward SDE and related $g$-expectation, \sl Pitman Res. Notes Math. \rm Ser. 364, Longman, Harlow, 1997, 141--159.

\bibitem{Peng-Shi 2000} S.~Peng and Y.~Shi, \it Infinite horizon forward-backward stochastic differential equations, \sl Stoh. Proc. Appl., \rm 85 (2000), 75--92.

\bibitem{Penrose 1955} R.~Penrose, \it A generalized inverse of matrices, \sl Proc. Camb. Phil. Soc., \rm 52 (1955), 17--19.

\bibitem{Sun-Yong 2018} J.~Sun and J.~Yong, \it Stochastic linear quadratic optimal control problems in infinite horizon, \sl Appl. Math. Optim., \rm 78 (2018), 145--183.

\bibitem{Sun-Yong 2020} J.~Sun and J.~Yong, \sl Stochastic Linear-Quadratic Optimal Control Theory: Open-Loop and Closed-Loop Solutions, \rm Springer, 2020.

\bibitem{Wei-Yu 2018} Q.~Wei and Z.~Yu, \it Time-inconsistent recursive zero-sum stochastic differential games, \sl Math. Control \& Related Fields, \rm 8 (2018), \rm

\bibitem{Yin 2008} J.~Yin, \it On solutions of a class of infinite horizon FBSDEs, \sl Stat.\& Probab. Lett., \rm 78 (2008), 2412--2419.

\bibitem{Yong 2006} J.~Yong, \it Remarks on some short rate term structure models, \sl J. Indust. Management Optim., \rm 2 (2006), 119--134.

\bibitem{Yong 2008} J.~Yong, \it A stochastic linear quadratic optimal control problem with generalized expectation, \sl Stoch. Anal. Appl., \rm 26 (2008), 1136--1160.

\bibitem{Yong-Zhou 1999} J.~Yong and X.~Y.~Zhou, \sl Stochastic controls: Hamiltonian systems and HJB equations, \rm Springer, 1999.

\bibitem{Yu 2017} Z.~Yu, \it Infinite horizon jump-diffusion forward-backward stochastic differential equations and their application to backward linear-quadratic problem, \sl ESAIM: Control, Optim. Cal. Var., \rm 23 (2017), 1331--1359.


\end{thebibliography}
\end{document}